\documentclass[12pt]{amsart} 
\usepackage[margin=1in]{geometry}
\usepackage{graphicx}
\usepackage{amsmath,amsfonts}
\usepackage{amssymb}
\usepackage{epstopdf}
\usepackage{caption,subcaption}
\usepackage{amsthm}
\usepackage{xcolor}
\usepackage{relsize}
\usepackage{float}
\usepackage[linesnumbered,lined,boxed]{algorithm2e}
\SetAlCapSkip{1em}
\usepackage{mathrsfs}
\usepackage{multicol}
\usepackage{enumitem}

\newtheorem{theorem}{Theorem}
\newtheorem*{theorem*}{Theorem}

\newtheorem{proposition}[theorem]{Proposition}
\newtheorem*{proposition*}{Proposition}

\newtheorem{remark}{Remark}[section]

\newcommand{\E}{\mathbb{E}}

\newcommand{\Lc}{\mathcal{L}}

\newcommand{\muES}{\mu_{\mathrm{ES}}}
\newcommand{\lambdaES}{\lambda_{\mathrm{ES}}}
\newcommand{\rhoES}{\rho_{\mathrm{ES}}}

\DeclareMathOperator*{\argmax}{argmax}

\title[Optimal Dividend Bands Revisited]{Optimal Dividend Bands Revisited: A gradient-based method and Evolutionary Algorithms}
\author[H. Albrecher]{Hansj\"org Albrecher} 
\address[H. Albrecher]{Department of Actuarial Science, Faculty of Business and Economics, University of Lausanne and Swiss Finance Institute, B\^atiment Extranef, Quartier UNIL-Dorigny, 1015 Lausanne, Switzerland}
\email{hansjoerg.albrecher@unil.ch}
\author[B. Garcia Flores]{Brandon Garcia Flores}\address[B. Garcia Flores]{Department of Actuarial Science, Faculty of Business and Economics, University of Lausanne, B\^atiment Extranef, Quartier UNIL-Dorigny, 1015 Lausanne, Switzerland}
\email{brandon.garciaflores@unil.ch}

\begin{document}
\bibliographystyle{abbrv}
	\maketitle
	\begin{abstract}\noindent 
		We reconsider the study of optimal dividend strategies in the Cramér-Lundberg risk model. It is well-known that the solution of the classical dividend problem is in general a band strategy. However, the numerical techniques for the identification of the optimal bands available in the literature are very hard to implement and explicit numerical results are known for very few cases only. In this paper we put a gradient-based method into place which allows to determine optimal bands in more general situations. In addition, we adapt an evolutionary algorithm to this dividend problem, which is not as fast, but applicable in considerable generality, and can serve for providing a competitive benchmark. We illustrate the proposed methods in concrete examples, reproducing earlier results in the literature as well as establishing new ones for claim size distributions that could not be studied before. 
	\end{abstract}
\section{Introduction}
Consider the optimal dividend problem for an insurance company whose surplus process evolves according to the Cramér-Lundberg model (see e.g.\ \cite{asmussen2010ruin}). The company pays dividends to shareholders in continuous time according to some admissible strategy $\pi$, and the objective is to identify the strategy that maximizes the expected sum of  discounted dividend payments until the event of ruin. If there are no constraints on the size of the payments, Gerber \cite{gerber1969entscheidungskriterien} showed that such a strategy always exists and is given by a \textit{band strategy} that partitions the interval $[0,\infty)$ into three sets $\mathscr{B}_0, \mathscr{B}_c$ and $\mathscr{B}_{\infty}$: whenever the current surplus level is in $\mathscr{B}_{0}$, no dividends are paid, when the current surplus level is in $\mathscr{B}_c$, all incoming premium is paid as dividends so that the same surplus level is maintained until the next claim arrives, and when the current surplus level is in $\mathscr{B}_{\infty}$, a lump sum payment to the first surplus level outside of $\mathscr{B}_{\infty}$ is carried out. This leads to a cascading strategy towards ruin (see e.g.\ the sample path illustration in Figure \ref{fig1} later).   Since Gerber's result, the optimal dividend problem and variants have been studied under many different and more general model assumptions and under various constraints (see for instance Avanzi \cite{avanzi09} and Albrecher \& Thonhauser \cite{at09} for surveys). \\

\noindent For many claim size distributions of practical relevance, the above band strategy turns out to collapse to a \textit{barrier strategy} (i.e., $|\mathscr{B}_c| = 1$), see e.g.\ Gerber \& Shiu \cite{gerber1998time}, and Avram et al.\ \cite{avram2007optimal} and Loeffen  \cite{loeffen2008optimality} for sufficient conditions on the Lévy measure under which a barrier strategy is optimal in the general framework of spectrally negative Lévy processes. In such a case the determination of the respective optimal barrier is rather straight-forward, if the scale function corresponding to the underlying surplus process is explicitly available (see e.g.\ Hubalek \& Kyprianou \cite{hubalek2011}). In contrast, already for the classical Cramér-Lundberg model it turns out to be surprisingly challenging to (even numerically) identify the sets $\mathscr{B}_0,\mathscr{B}_c,\mathscr{B}_\infty$ for more involved claim size distributions. \\

The value function of the optimal dividend problem has been identified as a viscosity solution of the corresponding Hamilton-Jacobi-Bellman equation by Azcue \& Muler \cite{azcue2005optimal}, which led to an iterative procedure for the numerical determination of the optimal bands (cf. \cite{azcue2005optimal,azcue2012optimal} as well as Schmidli \cite{schmidli2007stochastic} and Berdel \cite{berdel2014}). However, when following the respective algorithm, one faces two main difficulties:
\begin{enumerate}
\item\label{enum:Prob1} There is no complete understanding of the internal mechanism which furnishes the sets $\mathscr{B}_0, \mathscr{B}_c$ and $\mathscr{B}_\infty$ given an arbitrary set of parameters of the Cramér-Lundberg model, so that, as of now, concrete numerical solutions are known only for very few concrete and simple claim distributions. 
\item\label{enum:Prob2} The numerical approaches suggested in the literature require the solution of several differential equations, which can be very expensive in terms of computation time.
\end{enumerate}
While the numerical algorithms suggested so far were a by-product of the meticulous study of the existence and uniqueness of the solution of the optimal dividend problem and its characterization, in this paper we would like to take a different route. Relying on the fact that a band strategy is optimal and assuming that there are only finitely many such bands, we are interested to see if there are numerical alternatives to determine these optimal bands more efficiently and/or generally. In particular, we propose two respective numerical algorithms that differ from previous approaches.\\
The first one exploits the 'cascading' nature of band strategies and establishes a method based on gradients, when the value function is considered as a function of all band levels rather than the initial surplus level. This will lead to a rather fast numerical routine that can also be tailored to work for cases in which the scale function is not explicitly available, but its Laplace transform is (which is the case in general, since that Laplace transform is defined as the reciprocal of the (shifted) Laplace exponent of the underlying spectrally negative Lévy process). For every fixed number of bands, one then obtains an iterative procedure and can finally compare whether increasing the number of bands still improves the solution. \\
The second method also uses the explicit iterative formula for the expected discounted dividend payments for given band levels and uses it as the objective function in a general \textit{evolutionary algorithm}.  
Evolutionary strategies (ES) have been applied in the past with some success in reinsurance problems where the evaluation of the function to optimize is only possible through numerical procedures due to the non-existence of explicit algebraic expressions (see, for example, Salcedo-Sanz et al.\ \cite{salcedo2014effectively} and Roman et al.\ \cite{roman2018evolutionary}). In our context, we propose the use of an evolutionary self-adaptive strategy based on the algorithm originally proposed in the survey paper by Beyer \& Schwefel \cite{beyer2002evolution} which does not use the derivatives of involved functions and which can be easily implemented in common programming languages. This genetic algorithm is rather flexible and works particularly well in high-dimensional problems. We will hence adapt such an algorithm to the needs of our dividend problem, and indeed get to the same solutions as the other methods do. While in reasonably low dimension (like the dividend problem typically is, as there are only a few bands to consider) the computation time using this algorithm is not favorable when compared to the gradient-based approach, it is applicable in very general setups and can nicely serve as a benchmark for numerical procedures. Furthermore, it can also be useful in other application areas in risk theory, and since the idea and implementation of evolutionary algorithms may not be so commonly familiar in the risk theory community, we present the underlying principle and implementation variants in some detail in a separate section. \\

The rest of the paper is organized as follows. Section \ref{sec:2} contains some definitions and preliminaries on the model assumptions and the nature of the dividend band strategies. In Section \ref{sec:3} we summarize the previously available numerical procedures for the determination of optimal dividend bands. Section \ref{sec:4} then provides the expressions for the expected discounted payments that will be used in the numerical algorithms, in particular with respect to the 'cascading' view. Section \ref{sec:5} develops the gradient-based algorithm and discusses issues of its implementation. In Section \ref{sec:6} we give a general account of the idea behind evolutionary algorithms, and the necessary adaptations to the optimal dividend problem are discussed in Section \ref{sec:7}. Finally, in Section \ref{sec:num} we provide numerical illustrations. We first re-derive the known optimal bands for the well-known example of an Erlang(2) claim size distribution derived in Azcue \& Muler \cite{azcue2005optimal} as well as the mixture of Erlang claim size distribution established in Berdel \cite{berdel2014}. We then establish a new instance of a mixure of Erlang distributions for which a 4-band strategy is optimal. Subsequently we use our algorithms to derive the optimal barrier level in a risk model with Pareto claim sizes (where a barrier is known to be optimal due to Loeffen \cite{loeffen2008optimality}). We also implement an example with a mixture of Erlang and Pareto claims, which could not be handled with previously existing techniques and for which also two bands turn out to be optimal. Section \ref{sec:conc} concludes. 

\section{Definitions and preliminaries}\label{sec:2}
\begin{figure}[ht]
	\begin{center}
		\includegraphics[height=5cm]{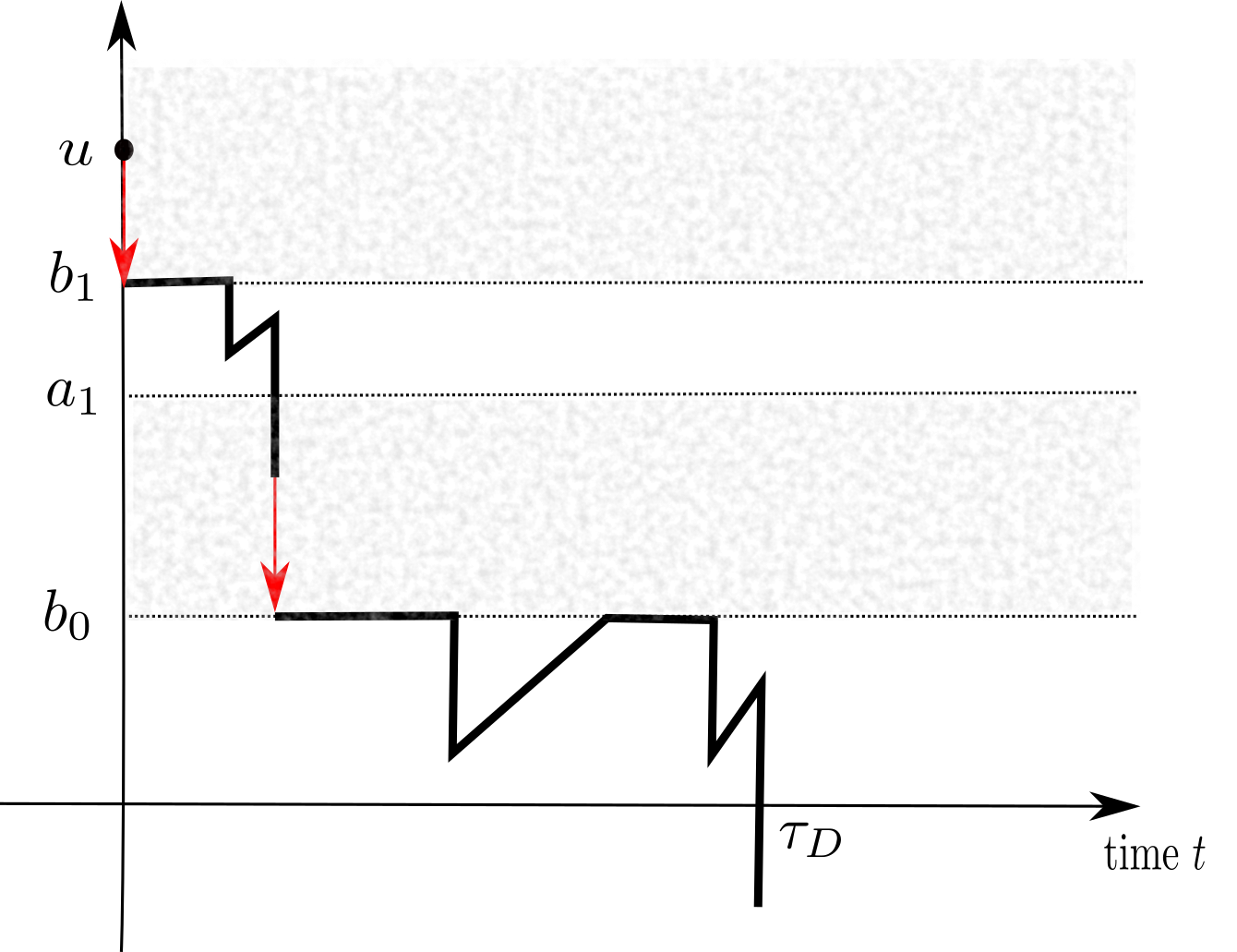}
	\end{center}
	\caption{A sample path with a 2-band strategy
	} \label{fig1}
\end{figure}

Consider the surplus process of an insurance portfolio in a Cramer-Lundberg model
\[
C_t = u + pt - \sum_{k=1}^{N_t}Y_k,\quad t\ge 0,
\]
with $(N_t)_{t\geq 0}$ a homogeneous Poisson process with rate $\lambda$ representing the arrival of claims, $(Y_k)_{k\geq 0}$ a sequence of absolutely continuous i.i.d.\ claim size random variables with density $f_Y$ and finite mean $\mu$, and the premium rate $p$ satisfying the positive safety loading condition $p=(1+\eta)\lambda \mu$ for some $\eta>0$. Let $(\mathscr{F}_t)_{t\geq 0}$ be the usual augmentation of the filtration generated by $(C_t)_{t\geq 0}$.\\
The dividend strategy $\pi$ is represented by the process $(U_t)_{\geq 0}$, where $U_t$ are the dividends paid up to time $t$. A dividend strategy is called \textit{admissible} if the associated process $(U_t)_{t\geq 0}$ is adapted to $(\mathscr{F}_t)_{t\geq 0}$, non-decreasing and càglàd. Denote by $\Pi$ the set of all admissible strategies.\\
For an admissible strategy $\pi\in\Pi$ we denote by $X_t = C_t-U_t$ the surplus after dividend payments. Let
\[
\tau_D = \inf\{t\geq 0 \mid X_t<0\}
\]
be the time until ruin, then
\[
V_\pi(u) = \E\left[\int _{0}^{\tau_D}e^{-\delta t}dU_t\mathlarger{\mathlarger{\mid}} X_0=u\right]
\]
is the expected value of the aggregated discounted dividends paid until ruin, where $\delta >0$ is the force of interest. The objective is then to identify the strategy that maximizes $V$ over all admissible strategies, that is, to find a strategy $\pi^\ast$ such that
\begin{equation}
V_{\pi^\ast} = \sup\limits_{\pi\in\Pi} V_{\pi}(x).\label{stochcontr}
\end{equation}
As pointed out in the introduction, the class of \textit{band strategies} is known to be optimal in this case (cf.\ \cite{gerber1969entscheidungskriterien}). A band strategy is defined by a partition of the positive half-line ${\mathbb R}=\mathscr{B}_{0}\cup \mathscr{B}_{c}\cup \mathscr{B}_{\infty}$ with the following properties:
\begin{itemize}
\item If $x\in \mathscr{B}_{0}$, there exists $\varepsilon>0$ such that $[x,x+\varepsilon)\subset \mathscr{B}_{0}$.
\item $\mathscr{B}_c$ is compact.
\item $\mathscr{B}_\infty$ is open in $[0,\infty)$.
\item If $x\not\in \mathscr{B}_{\infty}$ and there is a sequence $(x_{n})\subset \mathscr{B}_{\infty}$ such that $x_n\to x$, then $x\in \mathscr{B}_c$.
\item $(\sup\mathscr{B}_c,\infty)\subset \mathscr{B}_{\infty}$.
\end{itemize}
$\mathscr{B}_{0}$ corresponds to all surplus levels  for which $dU_t = 0$ (no dividends are being paid), $\mathscr{B}_c$ is the set of surplus levels for which $dU_t = p\,dt$ (all incoming premium is paid as dividends) and $\mathscr{B}_\infty$ is the set of surplus levels at which $U_{t+}-U_{t} = X_{t} - \sup\{b\in \mathscr{B}_{c}\mid X_t >b\}$ is applied (the smallest possible lump sum  is paid with which one leaves the set $\mathscr{B}_\infty$). As the focus of this paper is to provide alternatives to numerically identifying the optimal bands, we restrict ourselves to  finitely many bands, i.e., for given levels $a_0=0\leq b_0\leq a_1\leq \cdots \leq b_{m-1}$, the band strategy is given by
\[
\mathscr{B}_{0} = [0,b_0)\cup\bigcup_{k=1}^{m-1} [a_k,b_k), \; \mathscr{B}_c = \bigcup_{k=1}^{m-1} \{b_k\}, \; \mathscr{B}_{\infty} = \bigcup_{k=0}^{m-2} (b_k,a_{k+1})\cup (b_{m-1},\infty).
\]
We refer to this strategy as an \textit{m-band strategy} (see Figure \ref{fig1} for an illustration of a sample path with a 2-band strategy). \\
We conclude this section by making some remarks about notation: for a set $A\subset \mathbb{R}^{n}$, a function $f:A\to \mathbb{R}$, and a limit point $x\in A$, we denote by $f(x_1,\ldots,x_j-,\ldots,x_n)$ the limit $\lim_{y\to x} f(y)$ through points $y$ for which $y_j<x_j$, given that this limit exists and similarly for  $f(x_1,\ldots,x_j+,\ldots,x_n)$. The function $f$ is then continuous at $x$ if and only if $f(x_1,\ldots,x_j-,\ldots,x_n)=f(x_1,\ldots,x_j+,\ldots,x_n)$ for each $j=1,\ldots,n$, provided all right-hand and left-hand limits exist. Finally, in order to avoid cumbersome notation, we denote the partial derivative in the $i$-th variable by $D_{i}f$.

\section{The identification of bands in previous literature}\label{sec:3}
As mentioned in the introduction, the explicit identification of optimal bands has proved challenging even in the classical Cramér-Lundberg risk model. In the following we briefly summarize the available approaches in the literature. The typical approach is to derive the Hamilton-Jacobi-Bellman (HJB) equation 
\begin{equation}\label{HJB}
	\max\{1-V'(x),\Lc (V)(x)\} = 0,
\end{equation}
related to the  stochastic control problem \eqref{stochcontr}, where 
\begin{equation}\label{eq:diffgen}
	\Lc(f)(x) = pf'(x) - (\lambda+\delta)f(x) +\lambda\int_{0}^{x}f(x-y)f_{Y}(y)\;dy, \quad x>0
\end{equation}
is the infinitesimal generator of the discounted surplus process (see e.g.\ Azcue and Muler \cite{azcue2005optimal,azcue2012optimal}). Since $V$ is typically not sufficiently regular to satisfy the needs of \eqref{HJB} as a classical solution, one needs to look for viscosity solutions, and it can be shown that $V(x)$ is the unique viscosity solution of \eqref{HJB} satisfying a growth condition and a particular initial condition \cite{azcue2005optimal}. The numerical approach to find this solution is then an iterative procedure. In particular, Schmidli \cite{schmidli2007stochastic} proposed an algorithm for finding the levels of the optimal bands. This algorithm was formalized by Berdel \cite{berdel2014}, who considered the problem for the general case of phase-type claim distributions (cf.\ Algorithm 1).
\IncMargin{1em}
\begin{algorithm}[h]
\SetAlCapFnt{\sc}
\SetAlgoHangIndent{10em}

\SetKwInOut{Input}{Input}\SetKwInOut{Output}{Output}

\Input{Scale function $W_{\delta}$ and infinitesimal generator $\Lc$.}
\Output{Levels $B^*=(b_0^*,a_1^*,\ldots,b_{M-1}^*)$ of the best band strategy.}

\Begin{
$m:=0$\;
$b_0:=\sup\left\{x\geq 0 \mid W_{\delta}'(x) = \inf\limits_{y\geq 0} W_{\delta}'(y)\right\}$\;
$V_0(x) := \begin{cases}
  W_{\delta}(x)/W_{\delta}'(b_0)  & x\leq b_{0} \\
  V_{0}(b_0)+x-b_{0} & x>b_{0}
\end{cases}$\;
\While{$\Lc (V_{m})(x)>0$ for some $x>b_{m}$}{
$\mathscr{G}_{m}:=\{f^{a}_{m+1}, a>b_{m}\mid f^{a}_{m+1}(x)=V_{m}(x),\; x\leq a \text{ and}$
$\Lc (f)^{a}_{m+1}(x) = 0,\; x>a\}$\;\label{algline:Schmidli1}
$a_{m+1}:=\inf \left\{a> b_{m} \mid\inf\limits_{x>a} f^{a\prime}_{m+1}(x)=1\right\}$\;\label{algline:Schmidli2}
$b_{m+1}:=\sup \left\{x> a_{m+1} \mid f^{a_{m+1}\prime}_{m+1}(x)=1\right\}$\;\label{algline:Schmidli3}
$V_{m+1}(x) := \begin{cases}
  f^{a_{m+1}}_{m+1}(x) & x\leq b_{m+1} \\
  f^{a_{m+1}}_{m+1}(b_{m+1})+x-b_{m+1} & x>b_{m+1}
\end{cases}$\;
$m:=m+1$
}
}
\caption{Schmidli's algorithm}\label{alg:schmidli}
\end{algorithm}
\DecMargin{1em}
Here, $W_{\delta}$ is the scale function (see e.g.\ \cite[Ch.XI]{asmussen2010ruin}). 
When the optimal strategy is in fact a finite band strategy, Algorithm 1 is guaranteed to converge. However, depending on the particular distribution of the claims, the procedure can be computationally complex, as can be seen from Lines \ref{algline:Schmidli1} to \ref{algline:Schmidli3}. In each iteration of the algorithm, a family of functions $\mathscr{G}_{m}$ parametrized by the interval $(b_{m},\infty)$ is defined in such a way that for each $a>b_{m}$, the function $f^{a}_{m+1}$ solves \eqref{HJB} on $x>a$ with boundary condition $f^{a}_{m+1}(a) = V_{m}(a)$. Apart from some cases where this can be done explicitly, \eqref{HJB} has in general to be solved numerically. While this might not represent a problem for a couple of values of $a>b_{m}$, the difficulty arises when we consider Lines \ref{algline:Schmidli2} and \ref{algline:Schmidli3}, since they presuppose a full knowledge of the solutions in the entire interval $(b_{m},\infty)$ in order to compute the necessary extrema. As stated in \cite{berdel2014}, one can define $\bar{a}_{m+1}=\inf \left\{x> b_{m} \mid\Lc (V_{m})(x)>0\right\}$ and restrict the family to the interval $(b_{m},\bar{a}_{m+1})$. However, this introduces another extremum to be computed and one still has to consider the trade-off that arises at each step of the procedure when choosing a grid fine enough to discretize this new interval.\\

An alternative iterative algorithm for finding the optimal bands is proposed in Avram et al. \cite{avram2015gerber}, using stochastic sub- and super-solutions of \eqref{HJB} (their approach is formulated for general spectrally negative Lévy processes and the inclusion of fixed transaction costs with every dividend payment). Similarly to Algorithm \ref{alg:schmidli}, optimal levels are found in a sort of ``upwards'' approach finding higher band levels at each step of the procedure. However, instead of solving  integro-differential equations, each step consists of solving a stochastic control problem expressed through Gerber-Shiu functions. The advantage of this is that the problem is then reduced to finding the extrema of a low-dimensional function at each iteration, and there is no need for a full set of solutions of \eqref{HJB} as seen in, for example, Lines \ref{algline:Schmidli2} and \ref{algline:Schmidli3} of Algorithm \ref{alg:schmidli}. The Algorithm \ref{alg:alg_gradient} presented later in this paper sets out from a top-down approach, and then also leads to a bottom-up procedure that is formulated via discounted deficit densities explicitly (rather than general Gerber-Shiu functions), so that it eventually 
can be interpreted as a particular customization and implementation of the algorithm by Avram et al., see the details below.\\


\section{Properties of the Value Function}\label{sec:4}
In this section we recollect some properties of $V_{\pi}$, which will form the basis for the implementation of the numerical algorithms presented later.\\
For a fixed set of levels $a_0=0\leq b_0\leq a_1\leq \cdots \leq b_{m-1}$ of an $m$-band strategy $\pi$, we observe the following: for any $0\leq k\leq m-1$ and initial capital $u$ in $[a_{k},a_{k+1})$ (here $a_0=$ and $a_m=\infty$), the amount of dividends paid in a realization of the process will be the same as in a process with a barrier strategy with initial capital $u-a_{k}$ and barrier $b_{k}-a_{k}$, up until a claim makes the original process go below $a_k$. Denoting by $V_b$ the value function of a barrier strategy with barrier $b$, space-homogeneity and the Markov property imply that
\begin{equation}\label{eq:Eq1Vpi}
V_{\pi}(u) = V_{b_{k}-a_{k}}(u-a_{k})+ \mathbb{E}\left[V_{\pi}(a_k-D(u-a_{k},b_{k}-a_{k}))\right], \quad a_{k}\leq u< a_{k+1},
\end{equation}
where $D(u-a_{k},b_{k}-a_{k})$ denotes the deficit at ruin of a process with initial capital $u-a_{k}$, for which a barrier strategy with barrier $b_{k}-a_{k}$ is applied. In many cases, the density of the deficit at ruin can be computed by means of Gerber-Shiu functions (see \cite{gerber1998time,kyprianou2014fluctuations}) and the dividends-penalty identity (see \cite{lin2003classical,gly06}). Assume  henceforth that $D(u-a_{k},b_{k}-a_{k})$ has a density, which we denote by $f_D(\cdot,u-a_{k},b_{k}-a_{k})$. We observe that the variable inside the expectation in \eqref{eq:Eq1Vpi} is non-zero only when the deficit is at most $a_k$. We can therefore rewrite \eqref{eq:Eq1Vpi}, for any $k=0,\ldots,m-1$, as
\begin{equation}\label{eq:Eq2Vpi}
V_{\pi}(u) = V_{b_{k}-a_{k}}(u-a_{k})+ \int_{0}^{a_k}V_{\pi}(a_k-y)f_D(y,u-a_{k},b_{k}-a_{k})dy, \quad a_{k}\leq u< a_{k+1}.
\end{equation}
This set of equations provides the central formulas for computing the value of $V_\pi$ given a fixed set of levels: For $0\leq u< a_1$, the value of $V_\pi(u)$ is equal to $V_{b_0}(u)$, which is given in terms of the scale function of the process. We can then plug in these values in the integral in Equation \eqref{eq:Eq2Vpi} to obtain the value of $V_\pi(u)$ for $a_1\leq u <a_2$ and repeat this procedure in an iterative way to obtain the value of $V_\pi(u)$ for every $u$.
The problem of evaluating $V_{\pi}(u)$ is therefore reduced to the computation of the scale function $W_\delta$ and the density $f_D$. However, with knowledge of the scale function, the latter can be computed by means of the formula
\begin{equation}\label{eq:FDDef}
	\begin{split}
		f_D(y,u,b) = \lambda\int_{0}^{\infty}\left(\frac{W_{\delta}(u)W_{\delta}'(b-z)}{W_{\delta}'(b)}-W_{\delta}(u-z)\right)f_{Y}(y+z)\; dz.
	\end{split}
\end{equation}
see, e.g., \cite[Ch.X]{kyprianou2014fluctuations}. The setting in this last reference is that of general L\'evy processes. A more basic approach is to consider first the density $f_{D^0}(y,u)$ of the deficit at ruin with initial capital $u$ in the absence of a dividend strategy. Let $\Hat{f_Y}(s)$ denote the Laplace transform of the claim size density $f_Y$. Following e.g.\  \cite[Ch.XII]{asmussen2010ruin}, we know that $f_{D^0}(y,u)$ can be obtained as the inverse Laplace transform of 
\begin{equation}\label{eq:FSDef}
	\begin{split}
		\int_{0}^{\infty} e^{-su}f_{D^0}(y,u)\,du = \frac{\lambda(\Hat{w}(y,\rho)-\Hat{w}(y,s))}{ps-\delta-\lambda+\lambda \Hat{f_Y}(s)}, 
	\end{split}
\end{equation}
where $\Hat{w}:[0,\infty)^2\to \mathbb{R}$ is given by
\begin{equation}\label{eq:wLT}
\Hat{w}(y,s)= \int_{0}^{\infty} e^{-su}f_{Y}(y+u)\,du.
\end{equation}
From the dividends-penalty identity \cite{gly06} we then have
\begin{equation}\label{eq:DivPenIden}
f_D(y,u,b) = f_{D^0}(y,u)- \frac{W_\delta(u)}{W_{\delta}'(b)}D_2f_{D^0}(y,b).
\end{equation}
Note that we require $W_{\delta}$ to be differentiable in order to be able to use these formulas. Nonetheless, if for some $\alpha<1$ and $C>0$, we have $f_{Y}(x)\leq Cx^{-1-\alpha}$ for $x$ in some neighbourhood of the origin, then $W_{\delta}\in C^{q+2}(0,\infty)$ whenever $f_{Y}\in C^{q}(0,\infty)$ (see \cite{kuznetsov2012theory}), where $C^{q}(0,\infty)$ refers to the set of $q$-times continuously differentiable functions on the real positive line. A formula similar to \eqref{eq:FDDef} shows that the statement is also valid whenever we replace $W_{\delta}$ by $f_{D^0}$.\\
Equation \eqref{eq:Eq2Vpi} also reveals further properties of $V_{\pi}$ when we shift the focus from the initial capital $u$ to the band limits: for $m\geq 1$, we can identify the set of $m$-band strategies with the set
\[
\mathcal{B}_{m} = \{x\in \mathbb{R}^{2m-1}\mid 0\leq x_{1}\leq \cdots \leq x_{2m-1}\},
\]
and, for fixed $u>0$, we can consider the function $V^{m}:\mathcal{B}_{m}\to [0,\infty)$ given by $x\mapsto V_{x}(u)$. We have the following:
\begin{proposition}\label{prop:Prop}
If $f_Y \in C^{q}(0,\infty)$, $q\geq 2$, the function $V^{m}$ is continuously differentiable in the interior of the set 
\[
\mathcal{C}_{m}=\mathcal{B}_{m}\cap \{x\in \mathbb{R}^{2m-1}\mid  x_{2j-2}\neq u, j=1,\ldots,m\}.
\]
\end{proposition}
When $m=1$, we take $\{x\in \mathbb{R}^{2m-1}\mid  x_{2m-2}\neq u\}$ to be equal to $\mathbb{R}$.
\begin{proof}
We proceed by induction on $m$. For $m=1$, $\mathcal{C}_{1} = [0,\infty)$ and $\mathcal{D}_{1}=[0,u)\cup (u,\infty)$. For $0\leq u< b$, we have 
\[
V^{1}(b) = \frac{W_{\delta}(u)}{W_{\delta}'(b)}, \quad V^{1\prime}(b) = -\frac{W_{\delta}(u)W_{\delta}''(b)}{W_{\delta}'(b)^2},
\]
while for $0\leq b< u$,
\[
V^{1}(b) = u-b+V_{b}(b) = u-b+\frac{W_{\delta}(b)}{W_{\delta}'(b)}, \quad V^{1\prime}(b) = -\frac{W_{\delta}(b)W_{\delta}''(b)}{W_{\delta}'(b)^2}.
\]
Since $V^{1}(u-)=V^{1}(u+)=W_{\delta}(u)/W_{\delta}'(u)$ and $V^{1\prime}(u-)=V^{1\prime}(u+)=W_{\delta}(u)W_{\delta}''(u)/W_{\delta}'(u)^2$, the claim follows. Now, assume the claim is true for some $m\in\mathbb{N}$. We can write $\mathcal{C}_{m+1}$ as $\mathcal{C}_{m+1}=A\cup B\cup C$, where
\begin{align*}
A&=\{x\in \mathcal{C}_{m+1} \mid u<x_{2m-2}\},\\
B&=\{x\in \mathcal{C}_{m+1} \mid x_{2m-2}<u\leq x_{2m-1}\},\\
C&=\{x\in \mathcal{C}_{m+1} \mid x_{2m-1}<u\}.
\end{align*}
We observe the following: on $A$, $V^{m+1}=V^{m}\circ\pi$, where $\pi:\mathbb{R}^{2m+1}\to \mathbb{R}^{2m-1}$ is the projection onto the first $2m-1$ coordinates. Since $\pi$ maps the interior of $A$ into the interior of $\mathcal{C}_{m}$, by the induction hypothesis, it follows that $V^{m+1}$ is continuously differentiable in that set. Now, from \eqref{eq:Eq2Vpi} we have, for $b=(b_{0},a_{1},\ldots,b_{m})$ in $B$
\begin{equation*}
\begin{split}
V(b) = \frac{W_{\delta}(u-a_{m})}{W_{\delta}'(b_{m}-a_{m})}+ \int_{0}^{a_{m}}V^{m}_y(b)f_D(a_{m}-y,b_{m}-a_{m},b_{m}-a_{m})dy,
\end{split}
\end{equation*}
while for $b$ in $C$,
\begin{equation*}
\begin{split}
V(b) = u-b_{m}+\frac{W_{\delta}(b_{m}-a_{m})}{W_{\delta}'(b_{m}-a_{m})}+ \int_{0}^{a_{m}}V^{m}_y(b)f_D(a_{m}-y,b_{m}-a_{m},b_{m}-a_{m})dy.
\end{split}
\end{equation*}
From \eqref{eq:FSDef} and \eqref{eq:DivPenIden} we see that under the assumptions of the proposition, $f_D$ is twice continuously differentiable. Hence, using the induction hypothesis once again, it follows that $V$ is continuously differentiable in the interiors of $B$ and $C$. Moreover, since 
\begin{align*}
V(b_{0},a_{1},\ldots,u-)&=V(b_{0},a_{1},\ldots,u+),\\
D_{2m}V(b_{0},a_{1},\ldots,u-)&=D_{2m}V(b_{0},a_{1},\ldots,u+),
\end{align*}
and
\begin{equation*}
D_{2m+1}V(b_{0},a_{1},\ldots,u-)=D_{2m+1}V(b_{0},a_{1},\ldots,u+),
\end{equation*}
we conclude that $V$ is continuously differentiable in the interior of $\mathcal{C}_{m+1}$, concluding the proof.
\end{proof}
\begin{remark}\normalfont
By considering instead the set 
\[
\mathcal{D}_{m}=\mathcal{B}_{m}\cap \{x\in \mathbb{R}^{2m-1}\mid  x_{j}\neq u, j=1,\ldots,2m-1\},
\]
we can, in a similar manner, conclude that if $f_Y,W_{\delta}, f_{D^0} \in C^{q}(0,\infty)$, $q\geq 2$, then $V^{m}$ is $q-1$ times differentiable in the interior of $\mathcal{D}_{m}$, since in this case one does not have to consider the ``pasting'' at the points where $b_{j}=u$.
\end{remark}

\section{A gradient-based method}\label{sec:5}
From Proposition \ref{prop:Prop} and its proof, we can compute the gradient of the value function when we fix the initial capital $u$ and we look at it as a function of the levels. Assuming that $W_{\delta}$ and $f_{D^0}$ are twice differentiable and setting
\[
\mathcal{E}_m = \mathcal{B}_{m}\cap \{x\in \mathbb{R}^{2m-1}\mid x_{2m-2}<u\},
\]
we have for $b\in\mathcal{E}_{m}$ and $u<b_{m-1}$,
\begin{align}
	\begin{split}
		D_{2m-1}V^{m}(b) &= -\frac{W_{\delta}(u-a_{m-1})W_{\delta}''(b_{m-1}-a_{m-1})}{W_{\delta}'(b_{m-1}-a_{m-1})^{2}} \\
		&\quad + \int_{0}^{a_{m-1}}D_3f_D(a_{m-1}-y,b_{m-1}-a_{m-1},b_{m-1}-a_{m-1})V^{m-1}_y(b)dy,
	\end{split} \label{eq:PartialV2m1uLb}
	\\[0.5cm]
	\begin{split}
		D_{2m-2}V^{m}(b) &= \frac{W_{\delta}(u-a_{m-1})W_{\delta}''(b_{m-1}-a_{m-1})-W_{\delta}'(u-a_{m-1})W_{\delta}'(b_{m-1}-a_{m-1})}{W_{\delta}'(b_{m-1}-a_{m-1})^{2}}\\
		&\quad + \int_{0}^{a_{m-1}}D_1f_D(a_{m-1}-y,b_{m-1}-a_{m-1},b_{m-1}-a_{m-1})V^{m-1}_y(b)dy\\
		&\quad\quad - \sum_{i=2}^{3}\int_{0}^{a_{m-1}}D_if_D(a_{m-1}-y,b_{m-1}-a_{m-1},b_{m-1}-a_{m-1})V^{m-1}_y(b)dy\\
		&\quad\quad\quad + f_D(0,b_{m-1}-a_{m-1},b_{m-1}-a_{m-1})V^{m-1}_y(a_{m-1}),
	\end{split}\label{eq:PartialV2m2uLb}
\end{align}
while for $u>b_{m-1}$
\begin{align}
	\begin{split}
		D_{2m-1}V^{m}(b) &= -\frac{W_{\delta}(b_{m-1}-a_{m-1})W_{\delta}''(b_{m-1}-a_{m-1})}{W_{\delta}'(b_{m-1}-a_{m-1})^{2}} \\
		&\quad + \int_{0}^{a_{m-1}}D_3f_D(a_{m-1}-y,b_{m-1}-a_{m-1},b_{m-1}-a_{m-1})V^{m-1}_y(b)dy,
	\end{split} \label{eq:PartialV2m1uGb}
	\\[0.5cm]
	\begin{split}
		D_{2m-2}V^{m}(b) &= \frac{W_{\delta}(b_{m-1}-a_{m-1})W_{\delta}''(b_{m-1}-a_{m-1})}{W_{\delta}'(b_{m-1}-a_{m-1})^{2}} -1\\
		&\quad + \int_{0}^{a_{m-1}}D_1f_D(a_{m-1}-y,b_{m-1}-a_{m-1},b_{m-1}-a_{m-1})V^{m-1}_y(b)dy\\
		&\quad\quad - \sum_{i=2}^{3}\int_{0}^{a_{m-1}}D_if_D(a_{m-1}-y,b_{m-1}-a_{m-1},b_{m-1}-a_{m-1})V^{m-1}_y(b)dy\\
		&\quad\quad\quad + f_D(0,b_{m-1}-a_{m-1},b_{m-1}-a_{m-1})V^{m-1}_y(a_{m-1}),
	\end{split}\label{eq:PartialV2m2uGb}
\end{align}
and, in both cases, for $1\leq i < 2m-2$,
\begin{equation}\label{eq:PartialVi}
	D_iV^{m}(b) = \int_{0}^{a_{m-1}}f_D(a_{m-1}-y,b_{m-1}-a_{m-1},b_{m-1}-a_{m-1})D_iV^{m-1}_y(b)dy.
\end{equation}
We can solve these equations in an iterative manner to find candidate levels for the optimal band strategy whenever it is finite: call $b_{0}^{*}<a_{1}^{*}<\ldots <b_{M-1}^{*}$ the levels of the optimal strategy and assume for the moment that $b_{0}^{*}>0$. Since the first level has to occur at the largest global minimum of $W_{\delta}'$, we have $b_{0}^{*}=\sup\{x\geq 0 \mid W_{\delta}'(x) = \inf_{y\geq 0} W_{\delta}'(y)\}$. We observe that $W_{\delta}''(b_{0}^{*})=0$, so $D_{1}V^{1}(b_{0}^{*})=0$, independently of $u$.\\
If the barrier strategy at $b_{0}^{*}$ is not optimal and $u$ is such that $b_{1}^{*}<u<a_{2}^{*}$, $b_{0}^{*}$, $a_{1}^{*}$ and $b_{1}^{*}$ solve \eqref{eq:PartialV2m1uGb} to \eqref{eq:PartialVi} when $m=2$, since the global maximum is attained at the 2-band strategy with levels $b_{0}^{*}$, $a_{1}^{*}$ and $b_{1}^{*}$. Hence, since \eqref{eq:PartialVi} is zero regardless of $a_{1}$ and $b_{1}$ when $b_{0}=b_{0}^{*}$, we see that $a_{1}^{*}$ and $b_{1}^{*}$ are within the solutions to \eqref{eq:PartialV2m1uGb} and \eqref{eq:PartialV2m2uGb}. Moreover, these equations can be solved without regards to $u$, and, if we optimally choose the solution (so that we end up obtaining $a_{1}^{*}$ and $b_{1}^{*}$) we will see that $b_{0}^{*}$, $a_{1}^{*}$ and $b_{1}^{*}$ solve \eqref{eq:PartialV2m1uLb} to \eqref{eq:PartialVi} regardless of the value of $u>a_{1}^*$.\\
Continuing in this fashion, if a two-band strategy is not optimal, for $m=3$ and $i=2,3$, we have,
\begin{equation*}
	\begin{split}
		D_iV(b) &= \int_{0}^{a_{2}}f_D(a_{2}-y,b_{2}-a_{2},b_{2}-a_{2})D_iV^{2}_y(b)dy \\
		&= \int_{0}^{a_{1}^*}f_D(a_{2}-y,b_{2}-a_{2},b_{2}-a_{2})D_iV^{2}_y(b)dy\\
		&\quad+\int_{a_{1}^*}^{a_{2}}f_D(a_{2}-y,b_{2}-a_{2},b_{2}-a_{2})D_iV^{2}_y(b)dy.
	\end{split}
\end{equation*}
On the interval $(0,a_1^*)$, $V^{2}(b) = V^{1}(b)$ as functions of the initial capital, so $D_iV^{2}(b) = 0$. By the remarks of the previous paragraph, we also have $V^{2}_y(b)=0$ for all $y>a_{1}^*$, so we see that \eqref{eq:PartialVi} is always zero regardless of the value of $a_2$ and $b_2$. Hence, we can proceed again by solving \eqref{eq:PartialV2m1uGb} and \eqref{eq:PartialV2m2uGb}, choosing an 
optimal solution and test whether we proceed further with another band. At the $m+1$-th step, the equations that we need to solve can be written in a simpler form as
\begin{align}
	0 &= \frac{W_{\delta}(b_{m}-a_{m})W_{\delta}''(b_{m}-a_{m})}{W_{\delta}'(b_{m}-a_{m})^{2}}  - \int_{0}^{a_{m}}D_3f_D(a_{m}-y,b_{m}-a_{m},b_{m}-a_{m})V^{m}_{y}(b^*)dy, \label{eq:Grad1}
	\\
	1 &= f_D(0,b_{m}-a_{m},b_{m}-a_{m})V^{m}_{a_{m}}(b) + \int_{0}^{a_{m}}D_1f_D(a_{m}-y,b_{m}-a_{m},b_{m}-a_{m})V^{m}_{y}(b^*)dy \label{eq:Grad2}
\end{align}
where $b^*=(b_{0}^{*},a_{1}^{*},\ldots,b_{m-1}^{*})$.\\
If $b_{0}^{*}=0$, we can discard the equation for $D_{1}V^{m}$ and work instead on the interior of the set 
\[
\mathcal{E}'_{m}=\{x\in \mathbb{R}^{2m-2}\mid b_{0}^{*}\leq x_{1}\leq \cdots \leq x_{2m-3}\leq \min(u,x_{2m-2})\}
\]
by realizing that \eqref{eq:PartialV2m1uLb} to \eqref{eq:PartialVi} for $2\leq i<2m-2$ carry over verbatim to this situation. The same argument then shows that we can use the same procedure for obtaining the optimal levels. The procedure is described in Algorithm \ref{alg:alg_gradient}.

\IncMargin{1em}
\begin{algorithm}
	
	\SetAlgoHangIndent{5.6em}
	\SetAlCapFnt{\sc}
	
	\SetKwFunction{Solve}{solve}
	\SetKwFunction{Select}{select}
	\SetKwFunction{DefineV}{defineV}
	\SetKwFunction{Order}{order}
	\SetKwFunction{Selection}{selection}

	\SetKwInOut{Input}{Input}\SetKwInOut{Output}{Output}
	
	\Input{Scale function $W_{\delta}$ and density of deficit $f_{D}$.}
	\Output{Levels $B^*=(b_0^*,a_1^*,\ldots,b_{M-1}^*)$ of the best band strategy.}
	\BlankLine

	\Begin{
		
		$m:=0$\;
		$b_0:=\sup\left\{x\geq 0 \mid W_{\delta}'(x) = \inf\limits_{y\geq 0} W_{\delta}'(y)\right\}$\;
		$V^0(x) := \begin{cases}
			W_{\delta}(x)/W_{\delta}'(b_0)  & x\leq b_{0} \\
			V^{0}(b_0)+x-b_{0} & x>b_{0}
		\end{cases}$\;
		\While{$\Lc (V^{m})_x>0$ for some $x>b_{m}$}{
			$m:=m+1$\;
			$B^* =\, $\Solve{$D_{2m}V(b_0^*,a_1^*,\ldots,b_{m-1}^{*},a_{m},b_{m})=0,$
				$D_{2m+1}V(b_0^*,a_1^*,\ldots,b_{m-1}^{*},a_{m},b_{m})=0$}\;\label{algline:Alg3Solve2}
			$a_{m}^*,b_{m}^*:=\;$\Select{$B^*$}\;\label{algline:Alg3DefineAB}
			$V^{m}:=\,$\DefineV{$b_0^*,a_1^*,\ldots,b_{m}^{*}$}\;\label{algline:Alg3DefineV}
		}
	}
	\caption{Gradient-based algorithm for optimal dividends}\label{alg:alg_gradient}
\end{algorithm}\DecMargin{1em}

Algorithm \ref{alg:alg_gradient} starts by obtaining the first level $b_{0}^{*}$ of the optimal band strategy. If the barrier strategy at $b_{0}^{*}$ is optimal, the algorithm finishes. Otherwise, the algorithm enters into its main loop. After updating the number of bands, the loop proceeds to an application of an abstract \texttt{solve} function in Line \ref{algline:Alg3Solve2}, which is used to solve simultaneously Equations \eqref{eq:PartialV2m1uGb} and \eqref{eq:PartialV2m2uGb} (or equivalently, \eqref{eq:Grad1} and \eqref{eq:Grad2}) when the first $2m-1$ variables of $D_{2m}V$ and $D_{2m+1}V$ are fixed. In Line \ref{algline:Alg3DefineAB}, the best levels $a_{m-1}^{*}$ and $b_{m-1}^{*}$ are chosen by selecting the couple that produces the best value for the value function when the initial capital is set to $b_{m-1}^{*}$. Line \ref{algline:Alg3DefineV} uses the function \texttt{defineV} for creating the value function of the strategy with the levels $(b_0^*,a_1^*,\ldots,b_{m-1}^{*})$ found so far, finalizing the loop.\\
As explained before, Algorithm \ref{alg:alg_gradient} can be considered as a particular implementation of the algorithm proposed by Avram et al.\ \cite{avram2015gerber} but obtained after trying to solve the gradient equations in a sort of "backward" way. It is similar to Algorithm \ref{alg:schmidli} and Algorithm \ref{alg:alg_ods} in the sense that it is not possible to determine beforehand whether a finite band strategy is optimal or not. The advantage is, however, that one avoids having to fully specify the solutions to the HJB equation as in Algorithm \ref{alg:schmidli} while also avoiding the randomness involved in Algorithm \ref{alg:alg_ods}. 

\begin{remark}\normalfont
	Note that, in principle, one could derive explicit expressions that the optimal band levels should satisfy by means of equations \eqref{eq:Grad1} and \eqref{eq:Grad2}. However, even in the simple case of the claims following an Erlang(2) distribution, one already arrives at rather complicated expressions which involve combinations and products of exponentials and polynomials, and the resulting levels cannot be given in terms of elementary functions. However, these equations can still be solved through numerical methods, which is the basis of the gradient-based technique that is going to be introduced in the sequel.
\end{remark}

\section{Evolutionary Strategies}\label{sec:6}
Evolutionary strategies belong to a class of nature-inspired optimization algorithms which intend to mimic biological evolution by means of procedures roughly categorized as \textit{mutation}, \textit{recombination} and \textit{selection} procedures which incorporate tasks that resemble the way evolution is carried out in nature. Starting with a set of candidate solutions (called the \textit{parental population}), one produces a new set of candidate solutions by means of recombination and, through mutation, randomly alters it to form a second set of candidate solutions (called the \textit{offspring population}). One then uses selection to filter out the best candidate solutions from these two populations and iterates the process, replacing the previous parental population with the new population thus obtained. In general, recombination, mutation and selection tasks are problem-dependent and adjusted according to a diverse set of criteria. ES are classically referred by the way the offspring population is generated, and the notation for expressing it is the $(\muES/\rho_{ES}\overset{+}{,}\lambda_{ES})$-\textit{notation} (the subscripts ES are used here to differentiate these symbols from the previously defined $\mu$ and $\lambda$ in earlier sections). In this notation, $\muES$ refers to the size of the parental population at the beginning and end of each iteration, $\rho_{ES}$ refers to the amount of parents involved in the creation of one single offspring, randomly chosen without replacement, and $\lambdaES$ refers to the number of offsprings created in each iteration. The symbols ``+'' and ``,'' refer to the way selection is carried out: the first one indicates that the $\muES$ members of the new parental population are going to be extracted from a set obtained by merging the parental population together with the offspring population, while the symbol ``,'' indicates that the parental population is discarded after the creation of the $\lambdaES$ offsprings (so that in a $(\muES/\rho_{ES},\lambda_{ES})$ strategy we necessarily require $\lambdaES\geq \muES$).

The pseudo-code for the algorithm from Beyer  \cite{beyer2002evolution} is presented below in Algorithm \ref{alg:alg_beyer}, formulated in terms of a maximization problem. As suggested before, the basic objects handled by ES are populations, which in this strategy are modeled by tuples of the form $(x,s,f(x))$. In this representation, $x$ is simply a candidate solution belonging to the search space $X$. The element $s$ is used as a set of parameters aiding in the mutation procedure of the members of the population and leading the self-adaptive properties of the strategy. The last element is the value of the function to optimize at $x$, which needs to be stored in order to select elements in each iteration.

\IncMargin{1em}
\begin{algorithm}[h]
\SetAlCapFnt{\sc}

\SetKwFunction{Initialize}{initialize}
\SetKwFunction{Sample}{sample}
\SetKwFunction{SRec}{s\_recombination}
\SetKwFunction{XRec}{x\_recombination}
\SetKwFunction{SMut}{s\_mutation}
\SetKwFunction{XMut}{x\_mutation}
\SetKwFunction{Selection}{Selection}

\SetKwInOut{Input}{Input}\SetKwInOut{Output}{Output}

\Input{Function $f$ to maximize in unconstrained object space $E$.}
\Output{Solution to problem $x^\ast = \argmax\limits_{x\in E} f(x)$}
\BlankLine
\Begin{
$g:=0$\;
\Initialize{$\mathscr{P}^{(0)}:=\{(x_{0,k},s_{0,k},f(x_{0,k}))\mid k = 1,\ldots,\muES\}$}\;
\Repeat{$\mathrm{terminal\_condition}$}{
\For{$l:= 1$ \KwTo $\lambdaES$}{\label{algline:OffLoop}
$\mathscr{S}_l:=\,$\Sample{$\mathscr{P}^{(g)},\rhoES$}\;\label{algline:RecStart}
$\tilde{s}_l:=\,$\SRec{$\mathscr{S}_l$}\;
$\tilde{x}_l:=\,$\XRec{$\mathscr{S}_l$}\;\label{algline:RecEnd}
$\tilde{s}_l':=\,$\SMut{$\tilde{s}_l$}\;\label{algline:MutStart}
$\tilde{x}_l':=\,$\XMut{$\tilde{x}_l$,$\tilde{s}_l'$}\;\label{algline:MutEnd}
$\tilde{F}_{l}:=F(\tilde{x}_l')$\;\label{algline:ObjEva}
}
$\mathscr{O}^{(g)}:=\{(\tilde{x}_{l}',\tilde{s}_{l}',\tilde{F}_l)\mid l = 1,\ldots,\lambdaES\}$\;\label{algline:EndOff}
$\mathscr{P}^{(g+1)}:=\,$\Selection{$\mathscr{P}^{(g)}$,$\mathscr{O}^{(g)}$,$\muES$}\;\label{algline:Selection}
$g:=g+1$
}
}
\caption{The basic ES-algorithm}\label{alg:alg_beyer}
\end{algorithm}
\DecMargin{1em}

\begin{figure}[htb]
    \centering 
\begin{subfigure}{0.4\textwidth}
  \includegraphics[width=\linewidth]{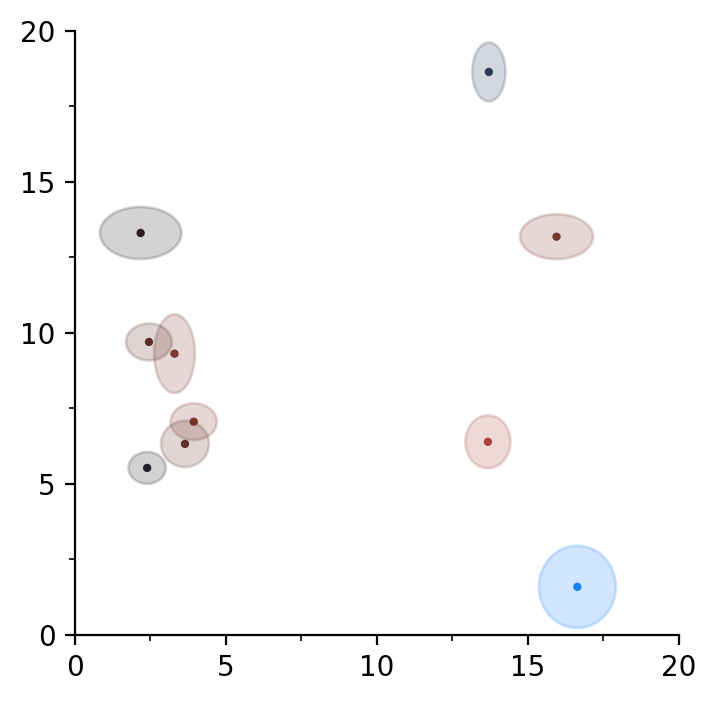}
  \caption{Initial Population}
  \label{fig:1}
\end{subfigure}\hfil 
\begin{subfigure}{0.4\textwidth}
  \includegraphics[width=\linewidth]{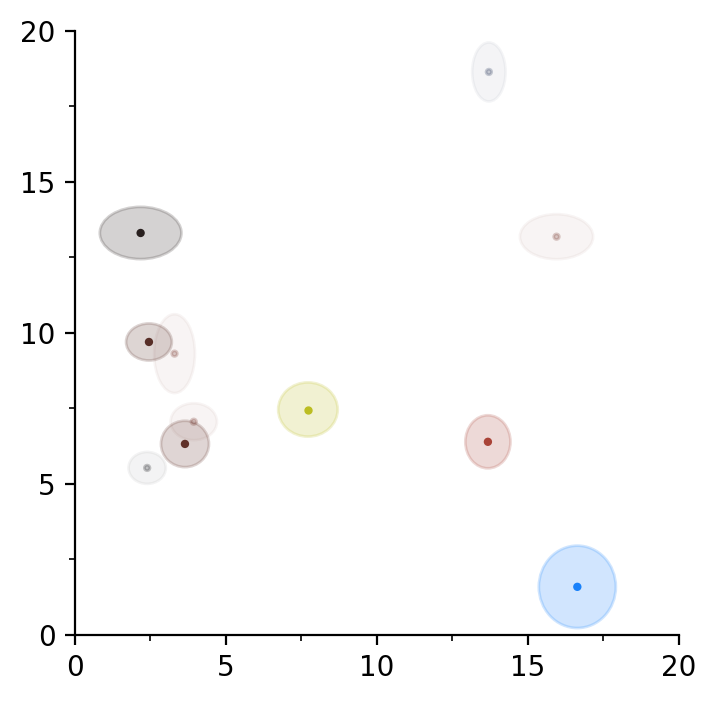}
  \caption{Recombination}
  \label{fig:2}
\end{subfigure}\hfil 

\medskip
\begin{subfigure}{0.4\textwidth}
  \includegraphics[width=\linewidth]{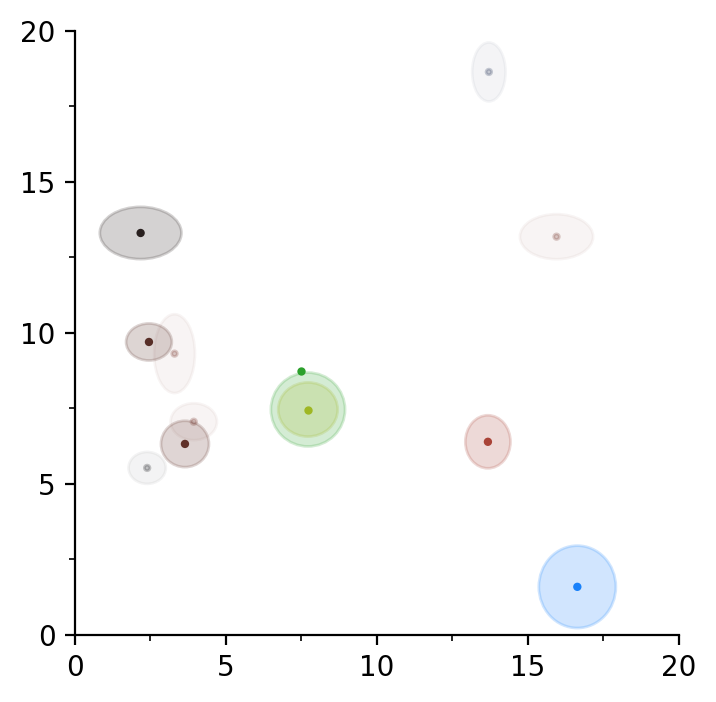}
  \caption{Mutation 1}
  \label{fig:3}
\end{subfigure}\hfil 
\begin{subfigure}{0.4\textwidth}
  \includegraphics[width=\linewidth]{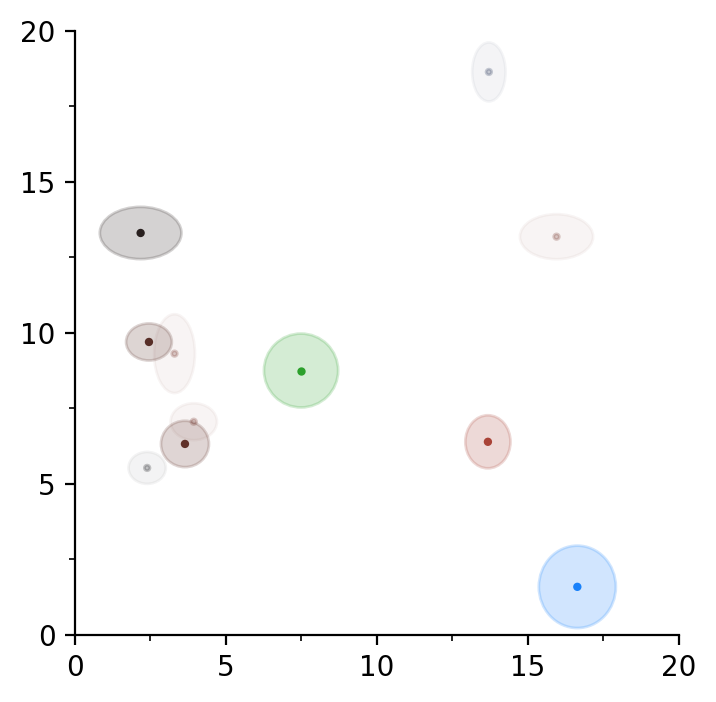}
  \caption{Mutation 2}
  \label{fig:4}
\end{subfigure}\hfil 

\medskip
\begin{subfigure}{0.4\textwidth}
  \includegraphics[width=\linewidth]{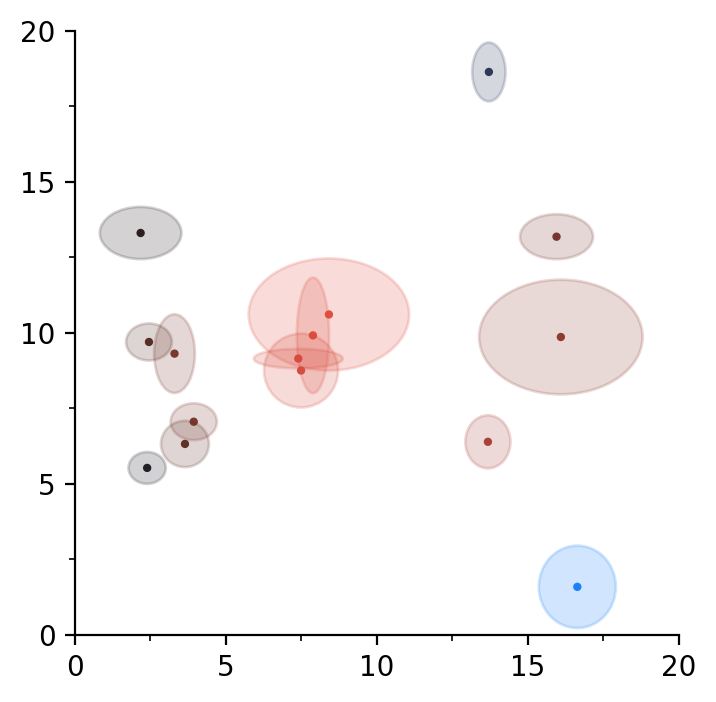}
  \caption{Evaluation}
  \label{fig:5}
\end{subfigure}\hfil 
\begin{subfigure}{0.4\textwidth}
  \includegraphics[width=\linewidth]{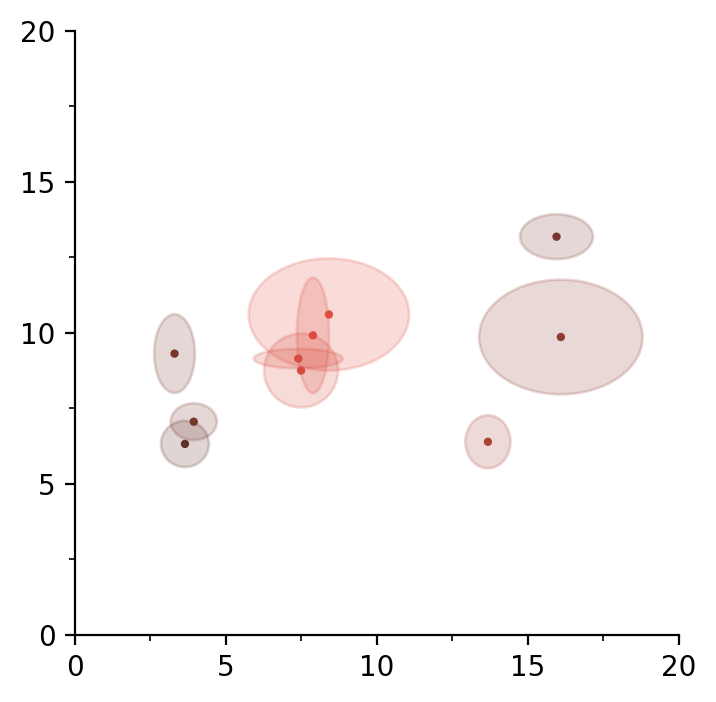}
  \caption{Selection}
  \label{fig:6}
\end{subfigure}\hfil 
\caption{Illustration of a $(10/5+5)$ ES.}
\label{fig:images}
\end{figure}

After initialization of the algorithm, which is usually carried out randomly, the algorithm enters into the main loop of the strategy for generating subsequent populations. This loop can roughly be described as an alternation of creating offsprings out of the parental population and selecting the replacing parental population out of these offsprings.\\
The process for creating offsprings is carried out from Line \ref{algline:OffLoop} to \ref{algline:EndOff}. Lines \ref{algline:RecStart} to \ref{algline:RecEnd} carry out the recombination procedure by first extracting a subsample from $\mathscr{P}^{(g)}$ of size $\rhoES$ without replacement. In real-valued spaces, a common recombination operator is the arithmetic mean, so, for example, $\texttt{s\_recombination}(\mathscr{S}_l) = \rhoES^{-1}\sum_{\mathscr{S}_l^s}s_{k}^{g}$ where $\mathscr{S}_l^s$ is the set of $s_k^{g}$ that belong to some tuple in $\mathscr{S}_l$. The mutation operator is then applied in lines \ref{algline:MutStart} to \ref{algline:MutEnd} by first mutating the strategy parameters and then the candidate solutions afterwards. While there is no established methodology for choosing the mutation operator, in \cite{beyer2001theory}, Beyer suggests that any operator should satisfy three requirements for a successful implementation of ES: \textit{scalability} (the ability to tune the strength of the mutation), \textit{reachability} (the ability to reach any other state $(x,s)$ within a finite number of steps) and \textit{unbiasedness}. Scalability is achieved by allowing the mutation of the object states, $x$, to be dependent on $s$. For $N$-dimensional real-valued search spaces, the parameter $s$ is generally used for controlling the variance of the mutation and in this regard, theoretical and practical considerations lead to a common mutation operator given by
\begin{equation}\label{eq:MutS}
\texttt{s\_mutation}(s_{l})^{j} = s_{l}^{j}\exp\left(\tau N_j\right), \quad j=1,\ldots,N,
\end{equation}
where $\tau$ is the learning-rate parameter and $N_j$ is a standard normally distributed random variable. Given the current parental state, the unbiasedness requirement simply means that the mutation procedure should not introduce any bias, and following the so-called \textit{maximum entropy principle}, this requirement immediately leads to mutation operators given by
\begin{equation}\label{eq:MutX}
\texttt{x\_mutation}(x_l,\tilde{s}_l)^{j}= x_l^{j}+\tilde{s}_{l}^{j}Z_j, \quad j=1,\ldots,N
\end{equation}
with $Z_j$ a standard normal random variable independent from the variables used to mutate $s$. However, Yao et al. suggest in \cite{yao1999evolutionary} and \cite{lee2004evolutionary} that, more generally, allowing $Z_j$ to have other kinds of stable distributions improves convergence speed and deals better with problems where several local extrema exist, dealing at once as well with a better handling of the reachability requirement.\\
The creation of the offspring is finalized in Line \ref{algline:ObjEva} by evaluating the objective function in the mutated objects $x$ and the offspring population is gathered in Line \ref{algline:EndOff}. The last step of the main loop is achieved in Line \ref{algline:Selection}, where the desired selection (plus or comma) takes place and the new parental population is created.\\
Figure \ref{fig:images} illustrates one iteration of a $(10/5+5)$ evolutionary strategy in a two-dimensional real space for the function $f(x,y)=\sqrt{\max\{10-(x-10)^2/2-(y-10)^2/2,0\}}$. Mutations for the parameters occur as in equations \eqref{eq:MutS} and \eqref{eq:MutX} with $N_j$ independent standard normal variables. Figure \ref{fig:1} represents the state of the population at the beginning of the iteration, where, in the notation of Algorithm \ref{alg:alg_beyer}, the $x_k^{g}$ are shown as the center of the ellipses, the $s_k^{g}$ as their axes and, using a blue-black-red scale, each point and ellipse is colored according to the value of $f$ at $x_k^{g}$. Figure \ref{fig:2} shows the first step in the creation of a single offspring: after randomly selecting 5 individuals from the original population (marked by the 5 darkest ellipses), the olive-colored ellipse is created after applying the recombination operator. In this case, recombination is given by the arithmetic mean, so that the olive point and the axes of the olive ellipse are the arithmetic means of the other 5 points and the ellipses' axes respectively. After recombination, mutation takes places, which is represented in Figure \ref{fig:3}. Here, the axes of the olive ellipse are mutated according to equation \eqref{eq:MutS}, generating the green ellipse. Sampling from the normal distribution centered in the olive point and variance given by the green ellipse, the green dot is created. Conclusion of the offspring individual's creation is depicted in Figure \ref{fig:4}, where the olive point and ellipse are deleted and the green ellipse is ``associated'' with the green point. After all offspring individuals have been created and their values according to $f$ have been computed, the parent and offspring population are merged, which is shown in Figure \ref{fig:5}. Finally, the \textit{plus} selection operator is used to discard the individuals with the lowest values of $f$, which finishes the iteration.

The presentation of the evolutionary strategy given so far raises two questions. First, in the case of real search spaces with the mutations defined in \eqref{eq:MutX}, Algorithm \ref{alg:alg_beyer} does not consider possible constraints imposed on the search space. A solution to this is the incorporation of restriction-handling techniques, like the inclusion of penalty functions, reparation of offspring, multiobjective optimization, etc., and the algorithm has to be adapted accordingly to the technique used (see \cite{kramer2010review} for an overview of several constraint-handling techniques in the context of ES).\\
The second question concerns the convergence of the algorithm. While theoretical results exist ensuring the almost sure convergence of the iterations (see, e.g., \cite{rudolph1996convergence}), the assumptions used in the statements of such results are usually quite restrictive or require a deep knowledge of the explicit form of the objective function, which leads to difficulties at the moment of the implementation. Despite this, evolutionary algorithms have been tested in a wide set of scenarios, proving to be effective tools for solving optimization problems.

\section{ES for the optimal dividend bands problem}\label{sec:7}
It is now of interest to see how competitive evolutionary strategies for the numerical determination of optimal band levels are in the present context.  Following the considerations from the previous section, in Section \ref{sec:num} we will use a $(\muES/\rhoES + \lambdaES)$ evolutionary algorithm to find the optimal band strategy for three distinct claim distributions: mixtures of Erlang distributions, a shifted-Pareto and a mixture of shifted-Pareto and Erlang distributions. As discussed earlier, only very few instances of explicit non-barrier optimal band strategies have been identified. In a what has become a classical example by now,   Azcue and Muler \cite{azcue2005optimal} identified a 2-band strategy for a case with Erlang(2,1) claims. Adding to this, in \cite{berdel2014} Berdel managed to expand this work by developing an algorithm for identifying non-barrier band strategies in the case of a mixture of Erlang distributions and some more general phase-type distributions. Our selection of mixtures of Erlang distributions for testing the ES was therefore made to compare its efficacy against an established baseline. Further, as can be seen from Section \ref{subsec:ObjectiveFuncEval}, the lack of explicit formulas for the scale function in the case of Pareto claims imposes the need of numerical approximations to the evaluation of $V_\pi$. As shown in Loeffen \cite{loeffen2008optimality}, for any choice of parameters, a barrier strategy is the optimal one for a Cramér-Lundberg model with shifted-Pareto claims. The second choice of claim distribution for the present work was then made to test the ES in a numerically-driven situation and test its respective efficacy. Finally, the mixture of Erlang and shifted-Pareto claim distribution was used as a means of testing the algorithm in uncharted territory.\\
Algorithm \ref{alg:alg_ods} displays the ES-algorithm adapted for the dividend-bands optimization, where $\odot$ stands for the element-wise multiplication operator, and $\bar{1}_k$ and $\bar{0}_k$ are the $k$-dimensional vectors of ones and zeros, respectively.\\
\IncMargin{1em}
\begin{algorithm}

\SetAlCapFnt{\sc}

\SetKwFunction{Initialize}{initialize}
\SetKwFunction{randn}{random\_normal}
\SetKwFunction{Sort}{sort}
\SetKwFunction{Order}{order}
\SetKwFunction{Selection}{selection}

\SetKwInOut{Input}{Input}\SetKwInOut{Output}{Output}

\Input{Initial capital $u_0$, upper bound on number of bands $M$, number of generations $G$ and variance bound $\varepsilon$ and value function $V_{\cdot}$.}
\Output{Levels $B^*=(b_0^*,a_1^*,\ldots,b_{M-1}^*)$ of the best $M$-band strategy.}
\BlankLine
\Begin{
\For{$m:= 1$ \KwTo $M$}{
$g:=0$\;
\Initialize{$\{B_{0,k}\mid k = 1,\ldots,\muES\}$}\;\label{algline:Alg2StartInit}
$s_{0,k}:=1.0,\; k=1,\ldots,\muES$\;
$\mathscr{P}^{(0)}:=\{(B_{0,k},s_{0,k},V_{B_{0,k}}(u_{0}))\mid k = 1,\ldots,\muES\}$\;\label{algline:Alg2EndInit}
\Repeat{$g=G$ or $\max(s_{g,0})<\varepsilon$}{
$\overline{s} = \frac{1}{\muES}\sum_{i=1}^{\muES}s_{g,i}$\;\label{algline:Alg2Recom1} 
$\overline{B} = \frac{1}{\muES}\sum_{i=1}^{\muES}B_{g,i}$\;\label{algline:Alg2Recom2}
$R:=\;$\randn{$1$}\;\label{algline:Alg2Random1}
\For{$l:= 1$ \KwTo $\lambdaES$}{
$s_R:=\;$\randn{$2m-1$}\;\label{algline:Alg2Random2}
$B_R:=\;$\randn{$2m-1$}\;\label{algline:Alg2Random3}
$\tilde{s}_l':=\,\overline{s}\odot \exp\left(s_R/\sqrt{4m-2}+R\cdot \bar{1}_{2m-1}/\sqrt{2\sqrt{4m-2}}\right)$\;\label{algline:Alg2StartMutation}
$\tilde{B}_l':=\,\max(\tilde{s}_l'\odot B_R+\overline{B},0_{2m-1})$\;
$\tilde{s}_l'':=\,$\Sort{$\tilde{s}_l'$,\Order{$\tilde{B}_l'$}}\;\label{algline:Alg2SortStart}
$\tilde{B}_l'':=\,$\Sort{$\tilde{B}_l'$}\;\label{algline:Alg2SortEnd}
$\tilde{V}_{l}:=V_{\tilde{B}_l''}(u_0)$\;\label{algline:Alg2EndRep}
}
$\mathscr{O}^{(g)}:=\{(\tilde{B}_{l}'',\tilde{s}_{l}'',\tilde{V}_l)\mid l = 1,\ldots,\lambdaES\}$\;
$\mathscr{P}^{(g+1)}:=\,$\Selection{$\mathscr{P}^{(g)}$,$\mathscr{O}^{(g)}$,$\muES$}\;\label{algline:Alg2Selection}
$g:=g+1$
}\label{algline:Alg2EndRepeat}
}
}
\caption{ES-algorithm for optimal dividends}\label{alg:alg_ods}
\end{algorithm}\DecMargin{1em}
Algorithm \ref{alg:alg_ods} is a $(\muES/\muES + \lambdaES)$-ES based on the basic strategy described in \cite{beyer2002evolution} for a search in a real unconstrained object space. Initialization is carried out in Lines \ref{algline:Alg2StartInit} to \ref{algline:Alg2EndInit}, where the function $\texttt{initialize}$ stands for random initialization of the candidate levels. Recombination of the parental population is done in Lines \ref{algline:Alg2Recom1} and \ref{algline:Alg2Recom2} using the arithmetic mean. In Lines \ref{algline:Alg2Random1}, \ref{algline:Alg2Random2} and  \ref{algline:Alg2Random3} the call of the function $\texttt{random\_normal}(k)$ represents the creation of an independent vector of dimension $k$ of standard normal random variables. Lines \ref{algline:Alg2StartMutation} to \ref{algline:Alg2SortEnd} show the implementation of the mutation operator, where each parent individual produces one offspring individual using log-normal multiplicative mutations for the exogenous parameters and normal mutations for the object parameters. The coefficients $1/\sqrt{4m-2}$ and $R\cdot \bar{1}_{2m-1}/\sqrt{2\sqrt{4m-2}}$ are learning rates which depend on the dimension of the search space and are based on both theoretical and empirical investigations. After the offspring has been created, repairing is carried out to ensure that the levels satisfy the condition $0\leq b_0\leq a_1<\cdots \leq b_{m-1}$. Line \ref{algline:Alg2SortStart} sorts the $s$ parameters according to the increasing order of $B$, while Line \ref{algline:Alg2SortEnd} sorts the object parameters in increasing order. Finally, the function in Line \ref{algline:Alg2Selection} performs \textit{plus} selection and outputs the population $\mathscr{P}^{(g+1)}:=\{(B_{g+1,k},s_{g+1,k},V_{B_{g+1,k}}(u_{0}))\mid k = 1,\ldots,\muES\}$ ordered in decreasing order according to the value of $V$, so when the terminal condition in Line \ref{algline:Alg2EndRepeat} is evaluated, $s_{g,0}$ holds the variances of the levels with the best fit.\\
Notice that the algorithm requires a value for the initial capital. While in principle this is a technical condition for the evaluation of $V_{\pi}$, caution should be taken: in case the optimal band strategy $\pi^\ast$ is finite with levels $b_{0}^{*}\leq a_{1}^{*}<\cdots <b_{m-1}^{*}$, for $b_{i}^{*}\leq u \leq a_{i+1}$, any other band strategy $\pi$ with first $i+1$ bands given by $b_{0}^{*}\leq a_{1}^{*}<\cdots <b_{i}^{*}$ and $u_0<a_{i+1}$ will satisfy $V_{\pi^*}(u_0)=V_{\pi}(u_0)$. Therefore, unless we can ensure $b_{m-1}^*\leq u_0$, any such $\pi$ will be the output of any optimization algorithm for which the initial capital is fixed. Following Lemma 3.3.1 in \cite{schmidli2006optimisation}, the inequality $b_{m-1}^*\leq u_0$ can be guaranteed by taking $u_0=p\lambda/(\delta(\lambda+\delta))$, which is the value that we use for all the iterations of the algorithm.\\
At this point, it is worthwhile to mention a key difference between the search method employed by this algorithm and the iterative algorithm discussed in, for example, \cite{azcue2005optimal} or \cite{azcue2012optimal}: given that the dimension of the search space has to be kept constant during the procedure, one has to fix in advance the number of bands for which the ES will try to identify the optimal levels. By observing that one can ``collapse'' levels in a band strategy, the $n$-bands strategies can be thought of as $m$-band strategies for $n\leq m$ and hence the algorithm would identify at once the best levels for all $n$-band strategies for $n\leq m$. If the optimal strategy is finite, one could then set $m$ large enough and use the ES to find the optimal levels. However, this approach requires the evaluation of $V_\pi$ for several bands and as explained in Section \ref{subsec:ObjectiveFuncEval} below, this is not efficient. Hence, a more efficient approach is instead to consecutively compute the the best 1-band, 2-band, 3-band, etc. strategies until collapsing of the levels is observed and then verify, by means of the HJB equation, that the proposed solution is in fact the optimal band strategy. Finally, efficiency is improved by skipping  the search for the optimal 1-band strategy and set $b_0:=\sup\left\{x\geq 0 \mid W_{\delta}'(x) = \inf_{y\geq 0} W_{\delta}'(y)\right\}$ in all searches.
\section{Numerical results}\label{sec:num}
We evaluate the performance of the procedures shown in the two previous sections by finding optimal band strategies for three study cases: claims distributed as a mixture of Erlang distributions, the case for a pure (shifted) Pareto distribution and a mixture between Erlang and Pareto distributions. The mixture of Erlang distributions is chosen because there are already explicit results available (see \cite{azcue2005optimal,berdel2014}) so that we can benchmark our algorithms. Given that no explicit expressions exist for the scale function when the claims follow a Pareto distribution, the second case was chosen to test the algorithms in a purely numerical situation (and in the case of a Pareto distribution, it is known that the optimal strategy is a barrier strategy, see \cite{loeffen2008optimality}). Finally, the mixture of Erlang and Pareto distributions was carried out to study the problem in a new context.
\subsection{Objective function evaluation}\label{subsec:ObjectiveFuncEval}
In the case where $f_Y$ comes from a mixture of Erlang distributions, the Laplace transform $\widehat{W_{\delta}}$ of the scale function is given in terms of a rational function, so explicit expressions in terms of the roots of the Lundberg equation can be found for $W_{\delta}$ and $f_D$. These expressions are then used for computing the value of $V_{\pi}$.\\
In the other two cases, numerical inversion of $\widehat{W_{\delta}}$ and \eqref{eq:FSDef} have to be carried out to find the values of $f_D$. Since the evaluation of this function is needed at several points, we opted for using a piece-wise linear approximation for $W_{\delta}$, $W_{\delta}'$, $f_{D^0}$ and its partial derivatives. The approximation was carried out in the following way: from the remarks of Section \ref{sec:7}, it can easily be seen that in the case where the optimal band strategy is finite, it is enough to restrict the domain of $f_D$ and $W_{\delta}$ to $[0,p\lambda/(\delta(\lambda+\delta))]^3$ and $[0,p\lambda/(\delta(\lambda+\delta))]$ respectively in order to find the optimal band levels. Hence, the approximation was done by evaluating 10,000 equidistant points in the interval $[0,p\lambda/(\delta(\lambda+\delta))]$ (including boundaries) and linear interpolation in between. The functions $f_{D^0}$ and $f_D$ were then computed using Equations \eqref{eq:FSDef} and \eqref{eq:DivPenIden}. Finally, except for the mixture of Erlang distributions appearing in Case I, the integrals appearing in \eqref{eq:Eq2Vpi}, \eqref{eq:PartialV2m1uGb}, \eqref{eq:PartialV2m2uGb} and \eqref{eq:PartialVi} were evaluated using numerical methods by using the \texttt{numpy}, \texttt{scipy} and \texttt{mpmath} libraries for Python 3. Equations \eqref{eq:Grad1} and \eqref{eq:Grad2} were solved by means of the MINPACK's hybrd and hybrj algorithms implemented in the  \texttt{scipy} library through the  \texttt{fsolve} function. \\
Finally, although there are no explicit results supporting the fact, the examples found in the literature indicate that the values for the different band values $b_1, b_2,\ldots$ are close to local minima of $W_{\delta}'$ found after the optimal $b_0$.\\
We now apply both numerical procedures introduced in this paper to the concrete examples. Below, the reporting times mean the clock time used to produce the results and do not include the time used for verification of the solution in an interval through the HJB equation, which is the same for both procedures.
\subsection{Case I: Erlang mixture claims}
The following three examples are considered (for which we know the explicit result already from \cite{azcue2005optimal} and \cite{berdel2014} for the first two):
\begin{itemize}
\item An $\mathrm{Erlang}(2,1)$ distribution with parameters $\lambda=10$, $\delta=0.1$ and $\eta=0.07$.
\item A mixture of the distributions $\mathrm{Erlang}(2,10)$, $\mathrm{Erlang}(3,1)$ and $\mathrm{Erlang}(4,0.1)$ with weights $0.025, 0.225$ and $0.75$ respectively, and parameters $\lambda=1$, $\delta=0.1$ and $\eta=0.405$.
\item A mixture of the distributions $\mathrm{Erlang}(2,10)$, $\mathrm{Erlang}(3,1.06775)$, $\mathrm{Erlang}(4,0.2325)$ and $\mathrm{Erlang}(5,0.05)$ with weights $0.005, 0.045$, $0.225$ and $0.725$ respectively, and parameters $\lambda=1$, $\delta=0.1$ and $\eta=0.4$.
\end{itemize}
For the first two distributions, a $(30/30 + 60)$-ES was used in both cases, with  bound on the variance equal to 0.01. For reasons that will be explained later, a combination of a $(30/30 + 60)$-ES and a $(1/1 + 1)$-ES was used for the third distribution, with same bound in the variance. The number of iterations vary from distribution to distribution.\\
For the Erlang(2,1) distribution, Table \ref{Table1} shows that we indeed find the optimal two-band strategy established in \cite{azcue2005optimal}.
\begin{table}[H]
\begin{tabular}{|c|c|c|c|c|c|}
\hline
  & Time & Iterations & $b_0$ & $a_1$ & $b_1$ \\ 
\hline
 Evolutionary Strategy & $\sim 1000$s & 1000 & 0 & 1.8064 & 10.2158\\  
\hline
 Gradient-Based & $<1$s & -- & 0 & 1.8030 & 10.2161\\
\hline
\end{tabular}
\caption{Results for the $\mathrm{Erlang}(2,1)$ distribution.}\label{Table1}
\end{table}
One can observe that the gradient-based approach is very fast, while the ES algorithm takes considerably longer time, but also arrives at the correct solution.\\
Table \ref{TableMixture1} shows the results for the first mixture of distributions, where a 3-band strategy is optimal. 
\begin{table}[H]
\begin{tabular}{|c|c|c|c|c|c|c|c|}
\hline
  & Time & Iterations & $b_0$ & $a_1$ & $b_1$ & $a_2$ & $b_2$  \\ 
\hline
 Evolutionary Strategy & $\sim 1$h  & $5000$ & $0.2617$ & $0.4668$ & $3.5249$ & $25.7390$ & $34.7857$\\  
\hline
 Gradient & $\sim 1$h & -- & $0.2615$ & $1.5230$ & $3.5246$ & $25.5763$ & $34.7696$\\
\hline
\end{tabular}
\caption{Results for the first mixture of Erlang distributions.}
\label{TableMixture1}
\end{table}
Several remarks for this example are in order. First, there is a clear discrepancy between the $a_1$-values obtained by the two methods. Comparing the value functions of both strategies shows that the strategy found by the ES provides higher values. However, the difference between the value functions is of the order of $10^{-4}$ while the norm of the gradient at that point is of the order $10^{-5}$ and the iterations do not reduce this value significantly, which shows why the gradient method has an early stop. Note that for this example, Berdel \cite{berdel2014} already studied the optimal dividend strategy, and her results are very similar to the ones in Table \ref{TableMixture1}, with only the values of $a_1$ and $a_2$ differing. The ES parameters above provide a larger value function, but again the difference is only of order $10^{-4}$. The explanation is that the step size in \cite{berdel2014} for solving the inf and sup in Lines 7 and 8 in Algorithm \ref{alg:schmidli} was set to be $10^{-4}$, whereas a smaller step size would have been needed to arrive at the above result.\\
As stated before, for the evolutionary strategy, the integrals in equations \eqref{eq:Eq2Vpi}, \eqref{eq:PartialV2m1uGb}, \eqref{eq:PartialV2m2uGb} and \eqref{eq:PartialVi} were not evaluated numerically but instead were computed exactly, by means of symbolical calculus in Mathematica. Now, for $k=2$, and $u\geq b_2$, \eqref{eq:Eq2Vpi} can be more explicitly written as 
\begin{align*}
	V_{\pi}(u) &= u - b_{2} + \frac{W_{\delta}(h_{2})}{W_{\delta}'(h_{2})}\\
	& + \int_{0}^{b_{0}} f_D(a_{2}-y,h_{2},h_{2})\frac{W_{\delta}(y)}{W_{\delta}'(b_{0})}\; dy\\
	& + \int_{b_{0}}^{a_{1}} f_D(a_{2}-y,h_{2},h_{2})\left(y-b_{0}+\frac{W_{\delta}(b_{0})}{W_{\delta}'(b_{0})}\right)\; dy\\
	& + \int_{a_{1}}^{b_{1}} f_D(a_{2}-y,h_{2},h_{2})\frac{W_{\delta}(y-a_{1})}{W_{\delta}'(h_{1})}\; dy\\
	& + \int_{a_{1}}^{b_{1}} f_D(a_{2}-y,h_{2},h_{2})\int_{0}^{b_{0}} f_D(a_{1}-z,y-a_{1},h_{1})\frac{W_{\delta}(z)}{W_{\delta}'(b_{0})}\; dz\; dy\\
	& + \int_{a_{1}}^{b_{1}} f_D(a_{2}-y,h_{2},h_{2})\int_{0}^{b_{0}} f_D(a_{1}-z,y-a_{1},h_{1})\left(y-b_{0}+\frac{W_{\delta}(b_{0})}{W_{\delta}'(b_{0})}\right)\; dz\; dy\\	
	& + \int_{b_{1}}^{a_{2}} f_D(a_{2}-y,h_{2},h_{2})\left(y-b_{1}+V_{\pi^{1}}(b_{1})\right)\; dy
\end{align*}
where $h_{i}=b_{i}-a_{i}$ and $V_{\pi^{1}}$ is the 2-band strategy obtained after deleting the last band from $V_{\pi}$. For mixtures of Erlang distributions, the scale function can be written as a linear combination of complex exponential functions with as many terms as roots of the Lundberg equation, assuming all of them are different. For the present case, there are 10 different roots. By means of formula \eqref{eq:FDDef}, it follows that $f_D(y,u,b)$ can be written as a sum of approximately 90 different terms involving $y$ with coefficients dependent on $u$ and $b$. Following this line of thought, the single integrals from the paragraph above have, in rough terms, 900 terms, while the first double integral has around 81,000 (in theory, further reductions that decrease these numbers considerably could in principle be possible, but the amount of terms implies that the human or computational effort for carrying out such operations is beyond reason). The computational effort for explicitly computing the integrals above was of around 1 hour, which implies that the total time for the ES was of around 2 hours, which doubled the computational time of the gradient method, but provided a slightly more accurate result. Although this procedure could also be carried out to test the computational time of the gradient method (which by virtue of the other cases would be expected to be smaller), we observe that the computational time for obtaining explicit expressions for \eqref{eq:PartialV2m1uGb} and \eqref{eq:PartialV2m2uGb} would be at least doubled, matching the current time of the ES.\\

\begin{figure}
  \begin{subfigure}{.5\linewidth}
  	\includegraphics[width=\linewidth]{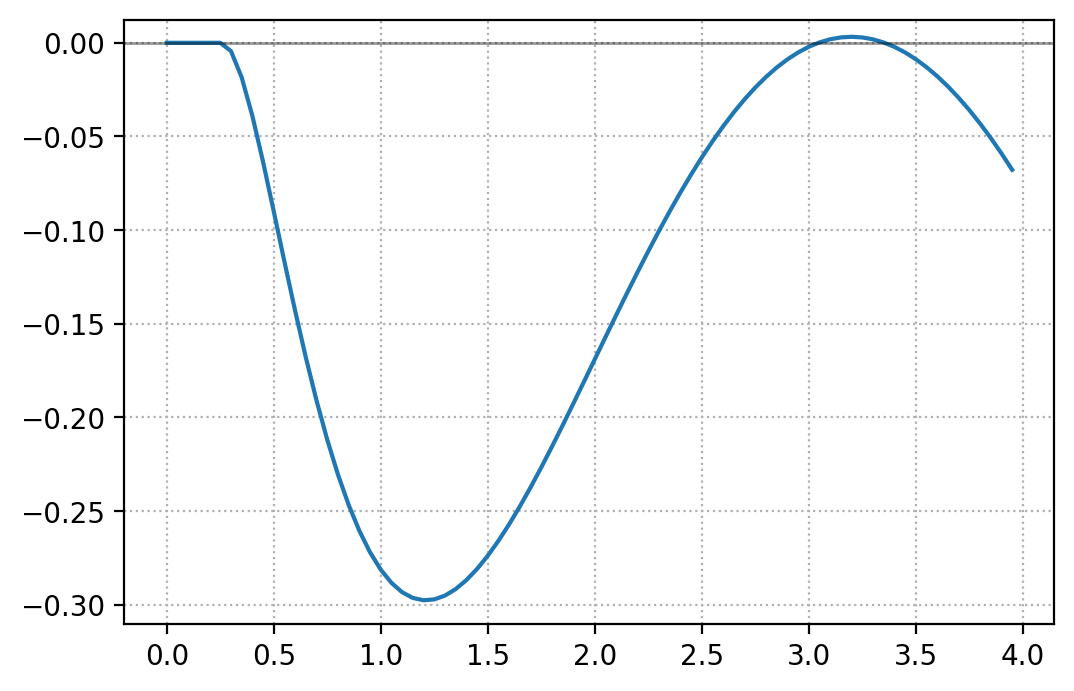}
  	\caption{HJB equation for the 1-band strategy}
  	\label{Subfig1}
  \end{subfigure}%
  \begin{subfigure}{.5\linewidth}
  	\includegraphics[width=\linewidth]{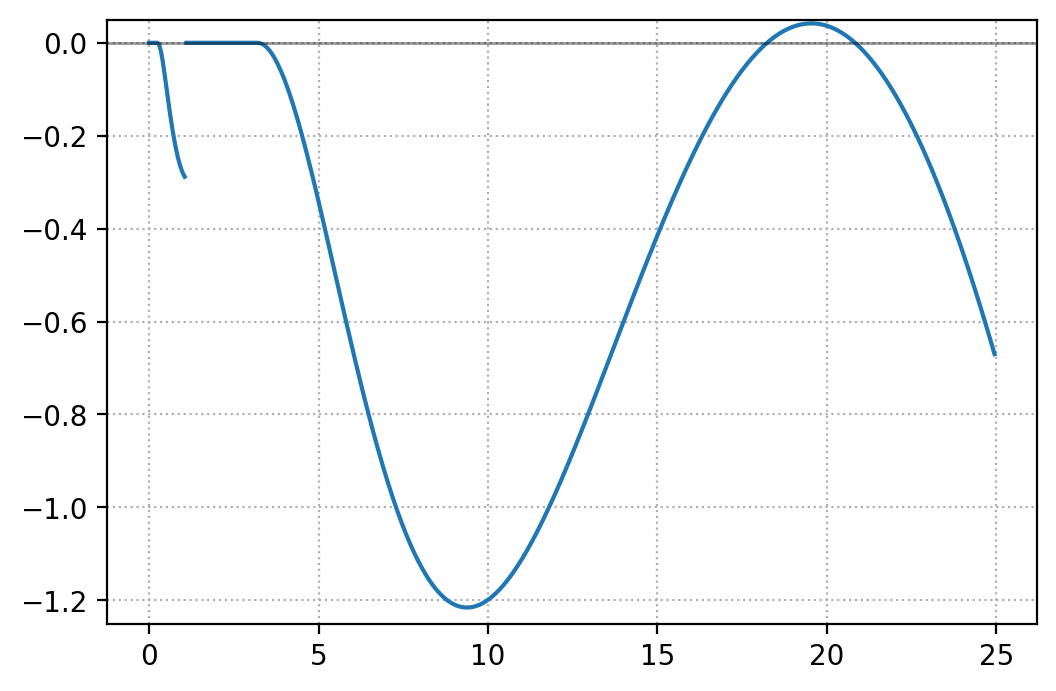}
  	\caption{HJB equation for the  2-band strategy}
  	\label{Subfig2}
  \end{subfigure}
  \begin{subfigure}{.5\linewidth}
    \vspace{0.5cm}
  	\centering
  	\includegraphics[width=\linewidth]{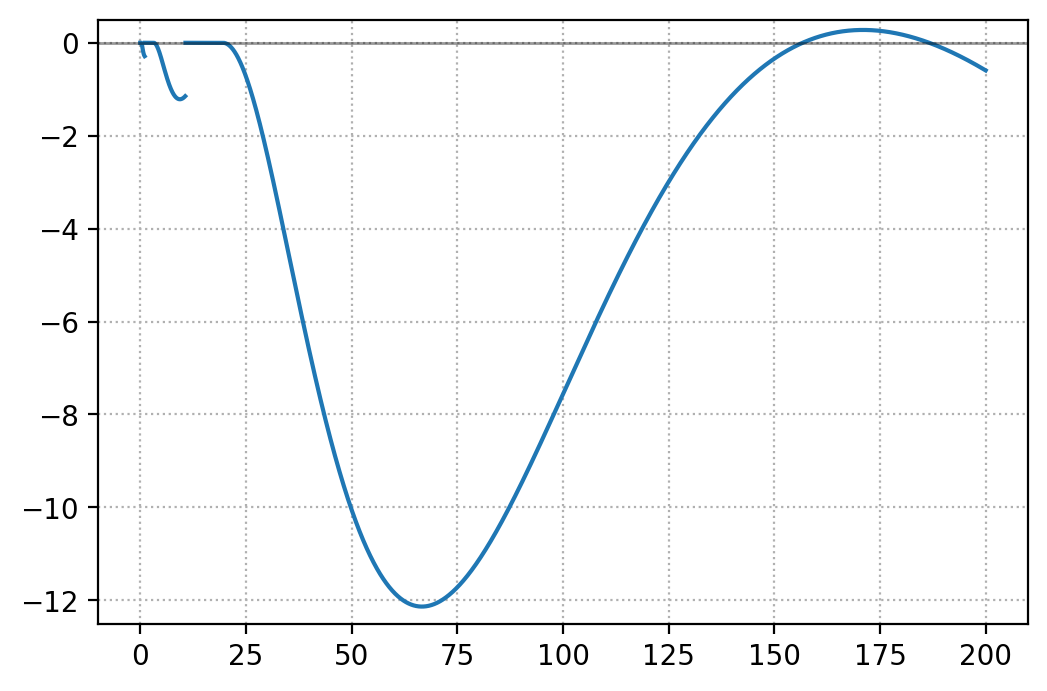}
  	\caption{HJB equation for the  3-band strategy}
  	\label{Subfig3}
  \end{subfigure}%
  \begin{subfigure}{.5\linewidth}
	\vspace{0.5cm}
  	\centering
  	\includegraphics[width=\linewidth]{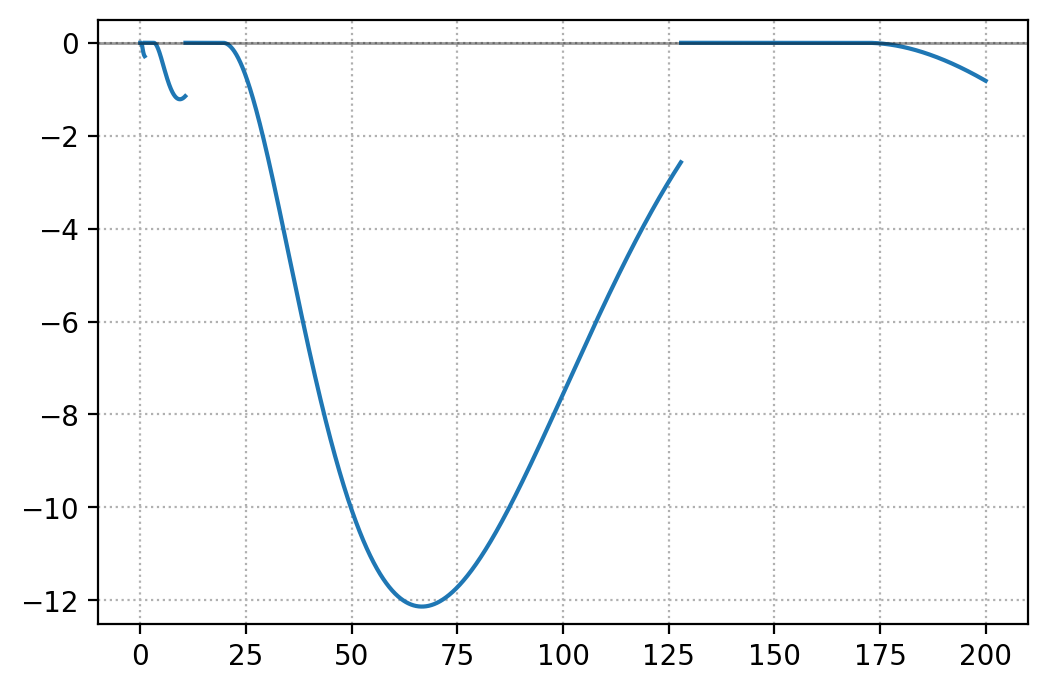}
  	\caption{HJB equation for the  4-band strategy}
  	\label{Subfig4}
  \end{subfigure}
  \begin{subfigure}{0.5\linewidth}
	\vspace{0.5cm}
  	\centering
  	\includegraphics[width=\linewidth]{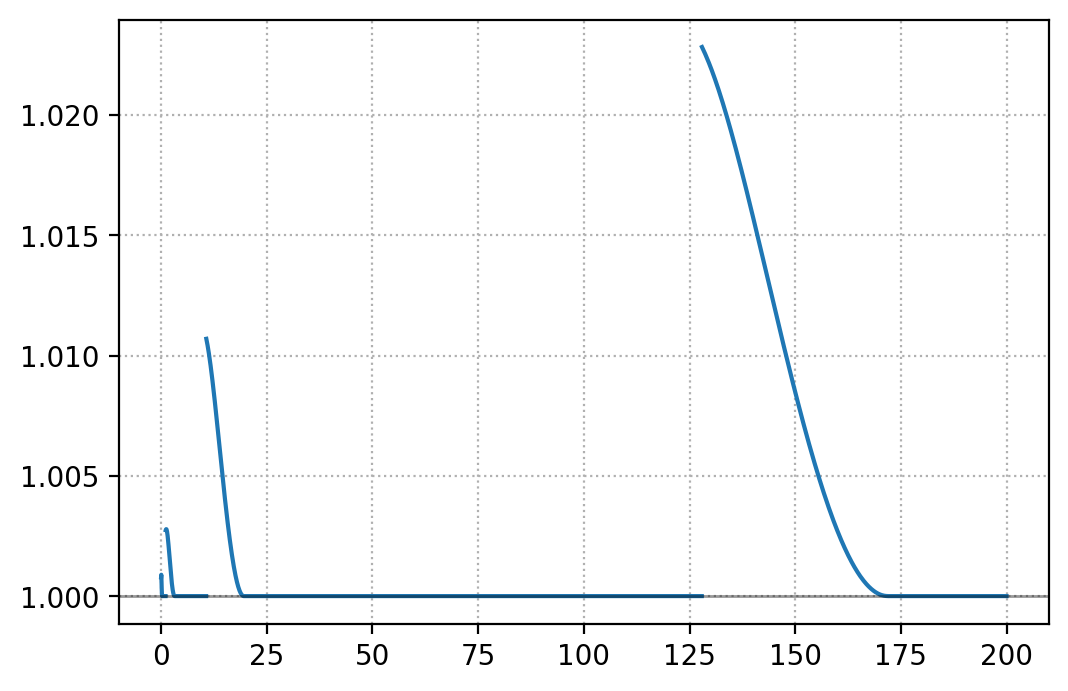}
  	\caption{Derivative of the value function}
  	\label{Subfig5}
  \end{subfigure}
  \caption{Plots of the l.h.s.\ of the HJB equation for the strategy with 1, 2, 3, and 4 bands based on the $a_i$'s and $b_i$'s from Table \ref{TableMixture2} as well as the derivative of the value function in the points where it exists.}
  \label{4BandHJB}
\end{figure}

Finally, using the intuition that the number of bands in the optimal strategy is related to the number of modes of the claim distribution, we were interested to establish a case where a 4-band strategy is optimal, and  the second mixture of Erlang distributions indeed leads to such an optimal 4-band strategy. The resulting optimal bands are given in Table \ref{TableMixture2}.
\begin{table}[H]
\begin{tabular}{|c|c|c|c|c|c|c|c|c|c|}
\hline
Time & Iterations & $b_0$ & $a_1$ & $b_1$ & $a_2$ & $b_2$ & $a_3$ & $b_3$  \\ 
\hline
$\sim 3$h  & $6000$ & $0.2562$, & $1.0543$ & $3.1988$ & $10.6647$ & $19.5499$ & $127.9288$ & $171.6044$\\  
\hline
\end{tabular}
\caption{Results for the second mixture of Erlang distributions.}
\label{TableMixture2}
\end{table}
These values were computed using only the ES technique. As before, the procedure was carried out in two steps, first using a $(30/30 + 60)$-ES for computing the values of the first three bands and consecutively using these values to reduce the problem to a two-dimensional optimization exercise, where a $(1/1 + 1)$-ES was used for computing the final two values. The reported time is for the combination of both procedures. The reason for proceeding in a two-step fashion was due to the observation that, as shown by these examples, the optimal values of of the $b_i$'s are usually located in the vicinity of the local minima of $W_{\delta}'$. The three smallest local minima are found in the interval $(0,20)$, which does not present any numerical complication. However, the last one is found at $172.7545$, which due to the nature of the scale function, produces exponentials with very large exponents at the moment of evaluating the value function, creating considerable numerical instabilities. To solve this issue, the value function and evolutionary strategy were re-implemented using arbitrary-precision floating-point-arithmetic, which decreased the speed at the moment of evaluating the value function. Since the $(1/1 + 1)$-ES is the strategy that requires the least evaluations of the objective function, it was chosen for obtaining the final values of $a_3$ and $b_3$. Figure \ref{4BandHJB} illustrates that this 4-band strategy is indeed the optimal strategy. Concretely, Figures \ref{Subfig1},\ref{Subfig2} and \ref{Subfig3} show that none of the first three strategies (with 1, 2 or 3 bands) is optimal, as the HJB equation attains positive values. Figure \ref{Subfig4} shows that this does not happen for the 4-band strategy and, moreover, Figure \ref{Subfig5} reveals that whenever the derivative exists, it is at least 1, so the solution is optimal.

\subsection{Case II: Pure Pareto claims}
Following \cite{loeffen2008optimality}, the optimal strategy will be a barrier strategy when claims have a shifted Pareto distribution with density function
\[
f_Y(y) = \alpha x_{0}^{-1}(1+x_{0}^{-1}y)^{-\alpha-1}, \quad y>0,
\]
and Laplace transform
\[
\Hat{f_Y}(s)= \alpha x_{0}^{\alpha}s^{\alpha}e^{sx_{0}}\Gamma(-\alpha,sx_{0}), \quad s>0,
\]
with $\Gamma$ the upper incomplete Gamma function and  $\alpha, x_{0}>0$. For the case at hand, we considered $x_0=1$ and $\alpha =1.5$, so that the claims have finite expectation and infinite variance. Moreover, the parameters of the Cramér-Lundberg process were taken to be $\lambda=10$, $\delta=0.1$ and $\eta=0.1$. The derivatives of the scale function were computed through their Laplace transforms and all of these were inverted using the de Hoog, Knight and Stones algorithm implemented in the library \texttt{mpmath}. The results are given in Table \ref{Table2}.\\
\begin{table}[H]
\begin{tabular}{|c|c|c|c|}
\hline
  & Time & Iterations & $b_0$  \\ 
\hline
 Evolutionary Strategy & $\sim 1000$s & 100 & 2.71036 \\  
\hline
 Gradient-Based & $<1$s & - & 2.71036\\
\hline
\end{tabular}
\caption{Results for the Pareto distribution.}\label{Table2}
\end{table}
Indeed, one arrives at an optimal barrier strategy, where for the evolutionary algorithm we only used 100 iterations to arrive at a running time that is comparable to the ones of the Erlang case, and the result is already well-aligned with the one of the gradient-based method. 
\subsection{Case III: Erlang and Pareto mixture claims}
Finally, let us consider a mixture of an $\mathrm{Erlang}(2,1)$ distribution and a shifted Pareto distribution ($\alpha=1.5$, $x_0=1$) with weights $0.8$ and $0.2$ respectively. The parameters of the Pareto distribution were chosen to match the mean of the Erlang component, while the weights were chosen to avoid a monotonicity of $W_{\delta}'$. The other parameters are again $\lambda=1$, $\delta=0.1$ and $\eta=0.1$. A $(150/150 + 100)$-ES was used and a 2-band strategy was found to be optimal. The results are shown in Table \ref{Table3}. 
\begin{table}[H]
\begin{tabular}{|c|c|c|c|c|c|}
\hline
  & Time & Iterations & $b_0$ & $a_1$ & $b_1$ \\ 
\hline
 Evolutionary Strategy & $\sim 8$h & 1000 & $0$ & $0.1524$ & $3.5115$\\  
\hline
 Gradient-Based & $\sim 1$h & - & $0$ & $0.0053$ & $3.8877$\\
\hline
\end{tabular}
\caption{Results for a mixture of an Erlang and a Pareto distribution.}\label{Table3}
\end{table}
In this case, the discrepancy between the results is more significant than in the other cases, with the gradient method providing a better solution. Due to the time it took for each evaluation, the ES was stopped after 1000 iterations, which meant that convergence was not fully achieved. The results from the gradient method also help to explain the poor performance of the ES: observe that the value of $a_1$ is rather close to zero, with an order of magnitude of $10^{-3}$. At the initialization of the algorithm, there is no knowledge of what the final order of magnitude will be, and during this experiment the initial values for the exogenous parameters were set to be 1, so, due to the projection into the zero for negative values, it takes many evaluations before the desired order of magnitude is achieved.
\section{Conclusions}\label{sec:conc}
In this paper we added two numerical alternatives to identify optimal dividend bands in the classical optimal dividend problem of risk theory. We illustrate that both of them are efficient, and their scope and applicability goes beyond the one of the previously discussed methods in the literature. The gradient-based method can be  particularly efficient. The second algorithm based on evolutionary strategies is satisfactory as well, and whereas in terms of computation times it can not compete with the gradient-based method for the complexity of this concrete problem, its range of applicability is even wider. In fact, ES algorithms can be an interesting competitor whenever an objective function can be efficiently evaluated, and it is known to work particularly well for higher-dimensional optimization problems, in which case the gradient alternative can be hard to explicitly compute or implement. We rederived optimal bands for some known cases, established new ones and also derived results for cases that were beyond the scope of previously available methods.  \\
The focus of this paper was on optimal dividend strategies in the Cramér-Lundberg model. However, since the equations used to derive the necessary functions for these two algorithms were obtained by means of Gerber-Shiu functions, one can in principle easily extend the range of applications to the case where a diffusion is added to the surplus process or even to the case where the surplus process is modelled by a generally spectrally-negative Lévy process satisfying the safety loading condition. Since evolutionary algorithms can be applied in rather general settings, it will also be interesting to see in future research other applications of this method in risk theory, particularly also in optimization problems with constraints, which may be handled with an introduction of a penalty term in the objective function (see e.g.\  \cite{kramer2010review}). 

\bibliography{Algorithm.bib}

\end{document}